\documentclass{amsart}
\usepackage{comment}
\usepackage{MnSymbol}
\usepackage{extarrows}
\usepackage[all,cmtip]{xy}
\usepackage{hyperref}

\DeclareMathOperator{\PGL}{PGL}

\DeclareMathOperator{\Cl}{Cl}
\DeclareMathOperator{\Clp}{Cl^+}
\DeclareMathOperator{\Disc}{Disc}
\DeclareMathOperator{\GL}{GL}

\DeclareMathOperator{\Gal}{Gal}

\DeclareMathOperator{\Norm}{Norm}

\DeclareMathOperator{\Tr}{Tr}
\DeclareMathOperator{\Frob}{Frob}
\DeclareMathOperator{\ord}{ord}

\newcommand{\Q}{{\mathbb Q}}
\newcommand{\Z}{{\mathbb Z}}
\newcommand{\C}{{\mathbb C}}
\newcommand{\R}{{\mathbb R}}
\newcommand{\F}{{\mathbb F}}

\newcommand{\cG}{\mathcal{G}}
\newcommand{\cF}{\mathcal{F}}
\newcommand{\cH}{\mathcal{H}}

\newcommand{\cI}{\mathcal{I}}

\newcommand{\cN}{\mathcal{N}}

\newcommand{\cJ}{{\mathcal J}}

\newcommand{\OO}{{\mathcal O}}
\newcommand{\ga}{\mathfrak{a}}

\newcommand{\sS}{\mathfrak{S}}

\newcommand{\ff}{\mathfrak{f}}

\newcommand{\fp}{\mathfrak{p}}
\newcommand{\fq}{\mathfrak{q}}
\newcommand{\mP}{\mathfrak{P}}

\begin {document}

\newtheorem{thm}{Theorem}
\newtheorem{lem}{Lemma}[section]
\newtheorem{prop}[lem]{Proposition}

\newtheorem{cor}[lem]{Corollary}
\newtheorem{conj}[lem]{Conjecture}
\newtheorem*{conj2}{Conjecture}

\theoremstyle{definition}

\theoremstyle{remark}

\title[]{Class field theory, Diophantine analysis and\\
the asymptotic Fermat's Last Theorem}

\author{Nuno Freitas}

\address{Mathematics Institute\\
	University of Warwick\\
	CV4 7AL \\
	United Kingdom}
\email{nunobfreitas@gmail.com}

\author{Alain Kraus}
\address{Sorbonne Universit\'e,
Institut de Math\'ematiques de Jussieu - Paris Rive Gauche,
UMR 7586 CNRS - Paris Diderot,
4 Place Jussieu, 75005 Paris, 
France}
\email{alain.kraus@imj-prg.fr}

\author{Samir Siksek}

\address{Mathematics Institute\\
	University of Warwick\\
	CV4 7AL \\
	United Kingdom}

\email{s.siksek@warwick.ac.uk}

\date{\today}
\thanks{Freitas is
supported by the
European Union's Horizon 2020 research and innovation programme under the Marie Sk\l{l}odowska-Curie grant 
agreement No.\ 747808 and the grant {\it Proyecto RSME-FBBVA $2015$ Jos\'e Luis Rubio de Francia}. 
Siksek is supported by an \emph{EPSRC LMF: L-Functions and Modular Forms Programme Grant EP/K034383/1}.}

\begin{abstract}
Recent results of Freitas, Kraus, \c{S}eng\"{u}n and Siksek,
give sufficient criteria for the asymptotic Fermat's 
Last Theorem to hold over a specific number field.
Those works in turn build on many deep theorems in arithmetic geometry.
In this paper we combine the aforementioned results 
with techniques from
class field theory, the theory of $p$-groups
and $p$-extensions, Diophantine
approximation and linear forms in logarithms,
to establish the asymptotic Fermat's Last Theorem
for many infinite families of number fields, and
for thousands of number fields of small degree.
For example, 
we prove the effective asymptotic Fermat's Last Theorem for the infinite family of fields 
$\Q(\zeta_{2^r})^+$ where $r \ge 2$.
\end{abstract}
\maketitle


\section{Introduction}

Let $K$ be a number field, and let $\OO_K$
be its ring of integers. 
The Fermat equation with prime exponent $p$ over $K$ is the equation
\begin{equation}\label{eqn:Fermat}
a^p+b^p+c^p=0, \qquad a,b,c\in \OO_K.
\end{equation}
A solution $(a,b,c)$ of \eqref{eqn:Fermat} is called trivial
if $abc=0$, otherwise non-trivial. 
The celebrated 
Fermat's Last Theorem, proved by Wiles \cite{Wiles}, asserts that the only
solutions to \eqref{eqn:Fermat} with $K=\Q$ and $p \ge 3$
are the trivial ones. The same statement, but for $p \ge 5$,
was proved for $\Q(\sqrt{2})$ by Jarvis and Meekin \cite{JarvisMeekin},
by Freitas and Siksek \cite{FS2} for a handful of other
real quadratic fields, and by Kraus \cite{Kraus} for the 
real cubic fields with discriminants $148$, $404$, $564$, 
and the quartic field $\Q(\zeta_{16})^+$.

The equation
\eqref{eqn:Fermat} defines a curve of gonality $p-1$,
and one can use this to show that there are non-trivial solutions over 
infinitely many number fields of degree $\le p-1$. 
Therefore a sharp
bound on~$p$ as in the works cited above does not exist for all number
fields. 
This remains true even if we restrict ourselves to totally real
number fields, as done for example in \cite{FS} and \cite{Kraus}. Indeed,
the field $\Q^{\mathrm{tr}}$, obtained by taking the 
union of all
totally real fields inside a given algebraic closure of $\Q$, 
is ample \cite{Pop}. Thus for
each
fixed exponent~$p$, the curve~\eqref{eqn:Fermat} has infinitely many points 
in~$\Q^{\mathrm{tr}}$. It thus becomes natural to consider the question asymptotically.

The \emph{asymptotic Fermat's Last Theorem over $K$} (or asymptotic FLT for
short)
is the statement that there is a bound $B_K$, depending only on the field $K$, 
such that for all prime exponents $p>B_K$, all solutions to~\eqref{eqn:Fermat} are trivial.
If $B_K$ is effectively computable,
we shall refer to this as the
\emph{effective asymptotic Fermat's Last Theorem over $K$}.
Let $\zeta_3$ be a primitive cube root of unity. The asymptotic
FLT is false if $\zeta_3 \in K$ as $(1,\zeta_3,\zeta_3^2)$
is a solution to \eqref{eqn:Fermat} for all $p \ne 3$. 
It seems reasonable to make the following
conjecture which is a consequence of the $abc$-conjecture for 
number fields (see \cite{browkin}).
\begin{conj2} 
Let $K$ be a number field. If $\zeta_3 \notin K$
then the asymptotic Fermat's Last Theorem holds for $K$.
\end{conj2}
In~\cite{FS}, Freitas and Siksek give a criterion for asymptotic FLT for totally real~$K$. This criterion can be formulated
in terms of the solutions of a certain $S$-unit equation, where
$S$ is the set of primes of $K$ above $2$, or equivalently
in terms of elliptic curves defined over $K$ having full $2$-torsion,
good reduction away from $S$ and with specified behaviour
at the primes in $S$. The proof 
builds on many deep results, including
modularity
lifting theorems over totally real fields due to Kisin, Gee and others,
Merel's uniform boundedness theorem, and Faltings' theorem for rational points
on curves of genus $\ge 2$. 
In~\cite{Haluk}, \c{S}eng\"{u}n and Siksek
establish a similar criterion for asymptotic FLT for
a general number field $K$,  subject to two standard
conjectures. In this paper we build on these results,
and with the help of class field theory, the theory of $p$-groups
and $p$-extensions, Diophantine approximation and linear
forms in logarithms, we 
establish asymptotic FLT for many infinite families of number fields,
and for thousands of number fields of small degree.
We remark that class field theory, and theory of cyclotomic fields,
was once considered the key to Fermat's Last Theorem \cite{LS}.
That approach was dramatically surpassed by the ideas 
of Frey, Serre, Ribet and Wiles. This paper demonstrates
that class field theory still has a r\^{o}le to play
in the beautiful story of Fermat.

\subsection*{A Conjecture of Kraus}
This paper is inspired by the following conjecture.
\begin{conj2}[Kraus~\cite{Kraus}]
Let $K$ be a totally real number field with narrow
class number~$1$. Suppose  $2$ is totally ramified in~$K$
and write $\mP$ for the unique prime above~$2$. Then
there are no elliptic curves over $K$
with full $2$-torsion and conductor~$\mP$.
\end{conj2}
\noindent Kraus showed that his conjecture 
implies asymptotic FLT for such $K$.
One objective  of this paper is to prove 
Kraus' conjecture in the following 
strengthened form.
\begin{thm}\label{thm:Krausgen}
Let $\ell$ be a rational prime. Let $K$ be a number field satisfying
the following conditions:
\begin{enumerate}
\item[(i)] $\Q(\zeta_\ell) \subseteq K$, where $\zeta_\ell$ is a primitive
$\ell$-th root of unity;
\item[(ii)] $K$ has a unique prime $\lambda$ above $\ell$;
\item[(iii)] $\gcd(h_K^+,\ell(\ell-1))=1$ where $h_K^+$ is the 
narrow class number of $K$.
\end{enumerate}
Then there is no elliptic curve $E/K$ with 
a $K$-rational $\ell$-isogeny, good reduction away from
$\lambda$ and potentially multiplicative reduction at $\lambda$.
\end{thm}
\noindent The proof can be found in Section~\ref{sec:Krausgen},
and makes use of class field theory, and the theory of $p$-groups
and $p$-extensions.

\subsection*{Asymptotic FLT for some infinite families of number fields}
The second objective of this paper is to use Theorem~\ref{thm:Krausgen}
to prove asymptotic
FLT for several infinite families of number fields.
Our results will be unconditional
for totally real fields. For number fields $K$
having at least one complex embedding, our theorems
are conditional on two standard conjectures.
We postpone the precise
statements of the two conjectures till Section~\ref{sec:criteria},
and now only briefly mention what they are.
\begin{itemize}
\item Conjecture~\ref{conj:Serre}: this is a weak version of
 Serre's modularity conjecture over general number fields. 
\item Conjecture~\ref{conj:ES}: this is a conjecture
in the Langlands Programme 
which says that every weight $2$ newform over $K$ with
rational integer Hecke eigenvalues has an associated elliptic
curve over $K$ or a fake elliptic curve over $K$.
\end{itemize}

\medskip

Before stating our main results, we need to introduce some notation. For a
number field $K$, we denote the class group by $\Cl(K)$
and the narrow class group by $\Clp(K)$. Their orders
are respectively the class number $h_K$, and the narrow
class number $h_K^+$. The class number $h_K$
is the degree of the Hilbert class field, which is the largest
abelian everywhere unramified extension of $K$.
The narrow class number $h_K^+$ is the degree of 
the narrow class field, 
which is the largest
abelian extension of $K$ unramified at all the finite places.
We shall
write $h^+_{K,2}$ for the largest power of $2$ dividing $h_K^+$.
This the degree of the narrow Hilbert $2$-class field of $K$,
which is the largest abelian $2$-extension of $K$
unramified at all the finite places. The degree $h^+_{K,2}$
is of course also the order of the $2$-Sylow subgroup of
$\Clp(K)$. For an ideal $\mP$ of $\OO_K$ we denote its
class in $\Clp(K)$ by $[\mP]$.

\begin{thm}\label{thm:main}
Let $K$ be a number field 
satisfying the following two hypotheses:
\begin{enumerate}
\item[(a)] $2$ is totally ramified in $K$;
\item[(b)] $h^+_{K,2}$ divides the order of $[\mP]$ in $\Clp(K)$,
where $\mP$ is the unique prime above~$2$.
\end{enumerate}
Then the following hold.
\begin{enumerate}
\item[(i)] If $K$ is totally real, then 
the asymptotic FLT holds over $K$.
Moreover, if all elliptic curves over $K$ with full $2$-torsion are modular, 
then
the effective asymptotic FLT holds over $K$.
\item[(ii)] If $K$ has at least one complex embedding, 
and we assume Conjectures~\ref{conj:Serre}
and~\ref{conj:ES} over~$K$ then the asymptotic FLT 
holds over $K$.
\end{enumerate}
\end{thm}
We can substantially strengthen the conclusion by adding one more
assumption to the theorem. By a $2$-extension $K^\prime$ of $K$
we mean a Galois extension of $K$ whose degree is a power of $2$.
\begin{thm}\label{thm:main2}
Let $K$ be a number field 
satisfying the following three hypotheses:
\begin{enumerate}
\item[(a)] $2$ is totally ramified in $K$;
\item[(b)] $h^+_{K,2}$ divides the order of $[\mP]$ in $\Clp(K)$,
where $\mP$ is the unique prime above $2$;
\item[(c)] $h_K$ is odd.
\end{enumerate}
Let $K^\prime/K$ be any $2$-extension unramified away from $\mP$
(and in particular unramified at the infinite places). 
\begin{enumerate}
\item[(i)] If $K^\prime$ is totally real, then 
the asymptotic FLT holds over $K^\prime$.
Moreover, if all elliptic curves over $K^\prime$ with full $2$-torsion are modular, 
then
the effective asymptotic FLT holds over $K^\prime$.
\item[(ii)] If $K^\prime$ has at least one complex embedding, 
and we assume Conjectures~\ref{conj:Serre}
and~\ref{conj:ES} 
over $K'$ then the asymptotic FLT holds over $K^\prime$.
\end{enumerate}
\end{thm}

We immediately obtain the following easy consequence of
Theorem~\ref{thm:main2}.  
\begin{cor} Let $K$ be a totally real number field
where $2$ is totally ramified. Suppose that $h_K^+$ is odd. 
Then asymptotic FLT holds over any totally real $2$-extension of $K$ unramified
away from~$2$.  
\end{cor}
\begin{proof}
Clearly, if the narrow class number $h^+_K$ is odd, then
$h_K=h_K^+$ and $h^+_{K,2}=1$ and so hypotheses (b), (c) are 
automatically satisfied. Hypothesis (a) is satisfied by assumption.
\end{proof}

\medskip

Let $r \ge 2$, and let $\zeta_{2^r}$ be a primitive $2^r$-th root
of unity. Write $\Q(\zeta_{2^r})^+$ for the maximal totally real
subfield of the cyclotomic field $\Q(\zeta_{2^r})$.
We now give a few consequences of Theorem~\ref{thm:main2}.
\begin{cor} \label{cor:oddhkp}
Let $K$ be a totally real number field satisfying (a), (b), (c).
Write $K_r$ for the compositum $K\cdot \Q(\zeta_{2^r})^+$. Then the asymptotic FLT holds for $K_r$.
\end{cor}
We note that the family $\{K_r\}$ are the layers in the 
$\Z_2$-cyclotomic extension of~$K$. Of course $K=\Q$
does satisfy (a), (b), (c). Here we can be more precise.
\begin{cor}\label{cor:cyclotomic}
The effective asymptotic Fermat's Last Theorem holds
over $\Q(\zeta_{2^r})^+$.
\end{cor}
\begin{proof}
For effectivity, we observe that $\Q(\zeta_{2^r})^+$
is contained in the $\Z_2$-extension of $\Q$,
and that modularity of elliptic 
curves over $\Z_p$-extensions of $\Q$ has been established
by Thorne \cite{Thorne}.
\end{proof}
Note that $K = \Q(\zeta_{2^r})^+$ is respectively $\Q$ and $\Q(\sqrt{2})$
for $r=2$, $3$.
In these
cases Corollary~\ref{cor:cyclotomic} is known in the stronger (non-asymptotic)
form
and is due respectively to Wiles \cite{Wiles} and to Jarvis and Meekin
\cite{JarvisMeekin}, as previously mentioned. 
For $r=4$ and $5$, it is proved by Kraus
\cite[Th\'{e}or\`{e}me 9]{Kraus} 
with $B_K=3$ and $B_K=6724$ respectively. 

\begin{cor}\label{cor:Lenstra}
Let $K/\Q$ be a $2$-extension unramified away from $2$ and $\infty$.
If $K$ has at least one complex embedding, assume 
Conjectures~\ref{conj:Serre} and~\ref{conj:ES}
over $K$. Then the asymptotic FLT holds for $K$.
\end{cor}
\begin{proof}
The fact that such $K$ satisfies conditions (a), (b)
is a result of 
Mark\u{s}a\u{\i}tis~\cite{Marksaitis} (we do however
reprove this as Corollary~\ref{cor:KrausNF}). 
\end{proof}

\begin{cor}\label{cor:correct}
Let $K$ be a real quadratic field.
Then $K$ satisfies conditions (a), (b), (c) of Theorem~\ref{thm:main2}
if and only if $K=\Q(\sqrt{2})$ or $K=\Q(\sqrt{\ell})$ or $K=\Q(\sqrt{2\ell})$
for some prime $\ell \equiv 3 \pmod{8}$. In particular,
for such $K$, the asymptotic FLT
holds over $K_r=K\cdot \Q(\zeta_{2^r})^+$ for all $r \ge 2$. 
\end{cor}
The proof of this corollary is given in Section~\ref{sec:quadratic}.
We remark in passing that the only imaginary quadratic
fields that satisfy conditions (a), (b), (c)
are $K=\Q(i)$ and $K=\Q(\sqrt{-2})$, as we explain
in Section~\ref{sec:quadratic}. Applying Theorem~\ref{thm:main}
to these two fields will merely give special cases of
Corollary~\ref{cor:Lenstra}, however  the stronger form of FLT over these
fields has already been established by
Turcas~\cite{Turcas}, subject only to Conjecture~\ref{conj:Serre}.

After the quadratic case, we 
establish asymptotic FLT for a very explicit
infinite family of totally real cubic fields.
This is done 
in Section~\ref{sec:cubic}.
\begin{thm}\label{thm:cubic}
Let $K$ be a totally real cubic field satisfying the following 
three conditions:
\begin{enumerate}
\item[(i)] $2$ is either totally ramified or inert in $K$;
\item[(ii)] $3$ ramifies in $K$;
\item[(iii)] the discriminant $\Delta_K$ is non-square (i.e.
the Galois group of the normal closure of $K$ is $S_3$).
\end{enumerate}
Then the asymptotic Fermat's Last Theorem holds for $K$.
\end{thm}

\begin{cor} \label{cor:cubic}
There is a positive
 proportion
of totally cubic real fields (ordered by discriminant)
satisfying asymptotic FLT.
\end{cor}
\begin{proof}
This now follows from \cite[Theorem 8]{BST} 
and \cite{Cohn}.
\end{proof}

The previous series of results 
establish asymptotic FLT
over infinitely many number fields. It is natural to wonder if we can 
establish asymptotic FLT over some number field of every
possible degree. We are able to give an affirmative though
conditional answer.
\begin{thm} 
Let $n \geq 2$ be an integer. There are infinitely many number fields of degree $n$ for which the asymptotic FLT holds subject to Conjectures~\ref{conj:Serre}~and~\ref{conj:ES}.

Moreover, if $n = 2^k$ with $k\ge 1$ or $n=3$ then
there are infinitely many totally real  fields of degree $n$ for which
asymptotic FLT holds (unconditionally).  
\end{thm}
\begin{proof}
For the second statement we take the fields $K_r$ in
Corollary~\ref{cor:correct} for varying~$\ell$, plus Corollary~\ref{cor:cubic}.
Thus we suppose $n\ge 5$.
By Theorem~\ref{thm:main}
all we need to do is show that there are infinitely many number
fields of degree $n$ such that $2$ is totally ramified
and having odd narrow class number. For this we need a special
case of a remarkable
theorem of Ho, Shankar and Varma \cite[Theorem 4]{HSV}.
Let $M/\Q_2$ be any totally ramified extension of degree $n$;
for example $M=\Q_2(\sqrt[n]{2})$. Let $r_2 \ge 1$, and $r_1=n-2r_2$.
The theorem of Ho, Shankar and Varma asserts the existence 
of infinitely many number fields $K$ of degree $n$, signature $(r_1,r_2)$,
and odd narrow class number, 
satisfying that $K \otimes \Q_2=M$. 

\end{proof}

\subsection*{A Computational Criterion for Asymptotic FLT}
The third objective of this paper is to provide a computationally viable
criterion for establishing  asymptotic FLT over specific number fields. 
The papers \cite{FS} and \cite{Haluk} give computational
criteria for asymptotic FLT in terms of the solutions
of a certain $S$-unit equation. There are algorithms
for determining the solutions to $S$-unit equations
(e.g.\ \cite[Chapter IX]{Smart}). However these 
algorithms require knowledge of the full unit group $\OO_K^\times$
of the number field $K$. Provably determining the full
unit group seems to be a computationally hard problem,
and is impractical in  current implementations 
(e.g.\ \texttt{Magma} \cite{magma}, \texttt{Pari/GP} \cite{PARI})
if the degree is much larger than $20$. It is however
much easier to determine a subgroup $V$ (say) of 
the unit group $\OO_K^\times$ of full rank
(\cite{Buchmann}, \cite{CDO}, \cite{Hanrot}).
Moreover, once one has a subgroup $V$ of full rank,
it is easy \cite[Section 5.3]{BF} for any given
prime $p$ to $p$-saturate $V$, i.e.\ to replace
$V$ by a larger subgroup of $\OO_K^\times$ whose index
is coprime to $p$. 
The following theorem gives a criterion
for asymptotic FLT that assumes knowledge not
of the full unit group but only
of a subgroup of full rank that is $2$-saturated.
\begin{thm}\label{thm:compCrit}
Let $K$ be a number field with one prime $\mP$
above $2$. Let $V$ be a subgroup of $\OO_K^\times$
such that the index $[\OO_K^\times : V]$ is finite and odd.
Let
\[
U:=\{ u \in V \; : \; u \equiv 1 \pmod{16\mP}\}.
\] 
Suppose every element of $U$ is the square of a unit.
Then the following hold.
\begin{enumerate}
\item[(A)] There is no elliptic curve $E/K$ with full $2$-torsion
and conductor $\mP$.
\item[(B)] If $K$ is totally real and $2$ is totally ramified
in $K$ then the asymptotic FLT holds over $K$.
\item[(C)] Suppose $K$ is a number field over which
Conjectures~\ref{conj:Serre} and~\ref{conj:ES} hold,
and $2$ is totally ramified in $K$. Then the
asymptotic FLT holds over~$K$. 
\end{enumerate}
\end{thm}
Theorem~\ref{thm:compCrit} is proved in Section~\ref{thm:compCrit}.
We remark that the criterion of the theorem is easy to test
computationally. One simply computes $U$ as the kernel
of the natural map $V \rightarrow (\OO_K/16\mP)^\times$,
and then tests whether each element in a chosen generating 
set is a square.
As an illustration, let
\begin{equation}\label{eqn:fn}
f_n(x)=\frac{1}{2\sqrt{-7}} \left( 
(1+\sqrt{-7})(x+\sqrt{-7})^n-(1-\sqrt{-7})(x-\sqrt{-7})^n
\right), \qquad n\ge 1
\end{equation}
The polynomial $f_n$ is monic, belongs to $\Z[x]$, and defines
a number field $K_n=\Q[x]/f_n(x)$ that is totally real and in
which $2$ totally ramifies (Lemma~\ref{lem:Kn}). Our computational
criterion establishes asymptotic Fermat over $K_n$ for
$1 \le n \le 6$, $8 \le n \le 14$, $15 \le n \le 20$, $23 \le n \le 27$,
$n=29$, $31$, $32$. This and other examples are found in 
Section~\ref{sec:examples}, where we also compare the relative
strength of Theorems~\ref{thm:main} and~\ref{thm:compCrit},
both computationally and theoretically.

\subsection*{Diophantine analysis and asymptotic FLT}
Whilst class field theory has distinguished historical
connections to Fermat's Last Theorem,
the subject of Diophantine analysis (Diophantine
approximation, linear forms in
logarithms) seems to have had
little or no influence on the mathematics
of the Fermat equation. In a surprising twist,
recent works (\cite{FS}, \cite{Haluk}) give
criteria for asymptotic FLT over certain
number fields conditional on properties
of the solutions of a specific $S$-unit equation.
Whilst the $S$-unit equation is treated by ad hoc methods
in \cite{FS}, \cite{Haluk}, and by class field theory for much of this
paper, the principal method of studying $S$-unit equations
is through Diophantine analysis (e.g. \cite{Smart}).
The final objective of this paper is to demonstrate 
that the methods of Diophantine analysis can be useful
in attacking asymptotic FLT over number fields.
\begin{thm}\label{thm:KrausconjII}
Let $\ell \equiv 1 \pmod{24}$ be a prime. 
The asymptotic FLT holds over~$\Q(\sqrt{\ell})$.
\end{thm}
The proof is given in Sections~\ref{sec:KrausconjII}--\ref{sec:dio}.
We remark that the family of quadratic
fields treated in this theorem is disjoint
from those treated in \cite{FS} or in the previous
parts of this paper. 
Indeed, in the family treated in Theorem~\ref{thm:KrausconjII},
the prime $2$ splits, 
whereas both in \cite{FS}
and in earlier parts of this paper, the focus
is on number fields with exactly one prime 
above $2$. 
That assumption is essential for the 
arguments of~\cite{FS}, and appears essential
to the class field theoretic arguments of this paper.

\subsection*{Acknowledgements}
We would like to thank Alex Bartel, Dominique Bernardi, Luis Dieulefait, Hendrik Lenstra, Bjorn Poonen, David Roberts and Haluk \c{S}eng\"{u}n for useful discussions.

\section{$p$-groups and $p$-extensions}
Let $p$ be a prime. A finite group $G$ is
said to be a $p$-group if its order $\#G$
is a power of $p$.
A finite extension
of fields $L/K$  is a $p$-extension if 
it is Galois and its degree $[L:K]$
is a power of $p$. Of course the Galois
group $\Gal(L/K)$ of a $p$-extension $L/K$
is a $p$-group. In this section we collect
well-known facts about $p$-groups and $p$-extensions that
we will make use of later.

The following is a standard result concerning $p$-groups;
see for example \cite[Section 2.1]{Suzuki}.
\begin{lem}\label{lem:maximal}
Let $G$ be a $p$-group. Then every maximal
subgroup of $G$ is normal of index $p$.
\end{lem}

We will denote the Frattini subgroup of a finite
group $G$ by $\Phi=\Phi(G)$; this is defined
as the intersection of all maximal subgroups of $G$.
The following is known as the \emph{Burnside basis theorem}
\cite[Theorem 2.1.16]{Suzuki}.
\begin{thm} \label{thm:suzuki}
Let $G$ be a $p$-group. Then $G/\Phi(G) \cong (\Z/p\Z)^r$
for some $r$, and so can be considered as an $r$-dimensional
$\F_p$-vector space. Let $x_1,\dotsc,x_s \in G$ 
and write $y_i=x_i \Phi(G) \in G/\Phi(G)$.
Then $G=\langle x_1,\dotsc,x_s \rangle$ if and only if 
$y_1,\dotsc,y_s$ span $G/\Phi(G)$. In particular, $G$
can be generated by $r$ elements.
\end{thm}

\begin{cor}\label{cor:abIsCyclic}
Let $G$ be a $p$-group, and write $G^\prime$ for its derived
subgroup. Suppose $G/G^\prime$ is cyclic. Then $G$ is cyclic
and $G^\prime=1$.
\end{cor}
\begin{proof}
By Theorem~\ref{thm:suzuki} we know that $G/\Phi(G)$ is abelian, and thus $\Phi(G) \supseteq G^\prime$.
We therefore have a natural surjection $G/G^\prime \rightarrow G/\Phi(G)$.
It follows that $G/\Phi(G)$ is cyclic. By Burnside's Basis Theorem
$G$ is cyclic.
\end{proof}

Finally we shall need the following standard result from Galois theory
for which we are unable to find a convenient reference.
\begin{lem}\label{lem:galClosure}
Let $L/K$ and $M/L$ be $p$-extensions. Let $N/K$
be the Galois closure of $M/K$. Then $N/K$ is a $p$-extension.
\end{lem}
\begin{proof}
Write $G=\Gal(N/K)$, $H=\Gal(N/L)$, $I=\Gal(N/M)$.
By the Galois correspondence and the hypotheses we have
\[
I \trianglelefteq H \trianglelefteq G,
\]
where the quotients $G/H$ and $H/I$ are $p$-groups. We
are required to show that $G$ is a $p$-group. 
Since $\# G = \#H \cdot \#\Gal(L/K)$
it is sufficient to show that $H$ is a $p$-group.

For $\sigma \in G$, we note that $\Gal(N/M^\sigma)=\sigma I \sigma^{-1}$. 
Observe, as $H$ is normal in $G$, 
that $\sigma I \sigma^{-1} \trianglelefteq H$ 
and the quotient $H/\sigma I \sigma^{-1}$ is isomorphic
to $H/I$, and so is a $p$-group.
As $N/K$ is the Galois closure of  $M/K$, we know that $N$ is generated
by the fields~$M^\sigma$, and so $\cap \sigma I \sigma^{-1}=1$. Therefore the natural
map
\[
H \rightarrow \prod_{\sigma \in G} H/\sigma I \sigma^{-1}
\]
is an injection. As the group on the right is a $p$-group, $H$ is also a $p$-group.
\end{proof}

Let $K$ be a number field. Recall that the Hilbert class field
is the largest abelian everywhere unramified extension of $K$;
its degree is the class number, which we denote by $h_K$. The narrow Hilbert
class field of $K$ is the maximal abelian extension of $K$ unramified away from
the infinite places; its degree is the narrow Hilbert class number which we
denote
by $h_K^+$.

We would like to thank Hendrik Lenstra for drawing our attention to the
following result. Part (a) is a result of Iwasawa \cite{Iwasawa}
(a proof is also found in \cite[Theorem 10.4]{Washington}). 
Part (b) is a straightforward
generalization.
\begin{thm}\label{thm:Iwasawa}
Let $K$ be a number field and let $\fq$
be a finite prime of $K$.
\begin{enumerate}
\item[(a)] 
Suppose $p \nmid h_K$.
Let $K^\prime/K$ be a $p$-extension unramified away from $\fq$
(and in particular, unramified at the infinite places).
Then $\fq$ is totally ramified in $K^\prime$ and $p \nmid h_{K^\prime}$.
\item[(b)] 
Suppose $p \nmid h_K^{+}$.
Let $K^\prime/K$ be a $p$-extension unramified away from $\fq$
and the infinite places.
Then $\fq$ is totally ramified in $K^\prime$ and $p \nmid h_{K^\prime}^+$.
\end{enumerate}
\end{thm}
\begin{proof}
We prove (b). 
The proof of (a) is almost identical. 
Write $G=\Gal(K^\prime/K)$.
Let $\fq^\prime$ be a prime of $K^\prime$ above $K$.
Let $I$ denote the inertia subgroup of $G$ for 
$\fq^\prime/\fq$. To show that $\fq$ is totally ramified
in $K^\prime$ is enough to show that $I=G$.
Suppose $I$ is a proper subgroup of $G$.
By Lemma~\ref{lem:maximal} there is a normal index $p$
subgroup $H$ of $G$ containing $I$. Consider
$K^{\prime\prime}={K^\prime}^H \subset K^\prime$. 
This is a Galois degree $p$ extension of $K$.
As $H$ contains $I$, the extension $K^{\prime\prime}/K$
 is unramified at some prime $\fq^{\prime\prime}$ above $\fq$
(and below $\fq^\prime$). But $K^{\prime\prime}/K$ is Galois,
and so it is unramified at all the primes above $\fq$.
It follows that $K^{\prime\prime}/K$ is a cyclic degree $p$ extension
unramified away from the infinite places, contradicting $p \nmid h_K^+$.
Therefore $\fq$ is totally ramified in $K^\prime$.

To complete the proof we would like to show that $p \nmid h_{K^\prime}^+$.
Suppose otherwise. So there is a cyclic degree $p$ extension $L/K^\prime$
unramified away from the infinite places. Now let $M/K$
be the Galois closure of $L/K$. This is a $p$-extension
by Lemma~\ref{lem:galClosure},
and it is unramified away from $\fq$ and the infinite places.
It follows from the first part that $\fq$ is totally ramified
in $M/K$. However, $M \supseteq L \supsetneq K^\prime \supseteq K$
and $L/K^\prime$ is unramified at any prime above $\fq$
giving a contradiction.
\end{proof}

The following is a result of Mark\u{s}a\u{\i}tis \cite{Marksaitis}.
It is 
immediate from part (b) of Theorem~\ref{thm:Iwasawa}.
\begin{cor}\label{cor:KrausNF}
Let $K/\Q$ be a $2$-extension unramified away from $2$, $\infty$.
Then $K$ has odd narrow class number and  $2$ totally 
ramifies in $K$.
\end{cor}

\section{Proof of Theorem~\ref{thm:Krausgen}}\label{sec:Krausgen}
Let $\ell$ be a rational prime.
Let $K$ be a number field satisfying conditions (i)--(iii)
of Theorem~\ref{thm:Krausgen}. 
In particular, there is
a unique prime $\lambda$ of $K$ above $\ell$. 
Let $G_K=\Gal(\overline{K}/K)$ 
be and 
$I_\lambda \subset G_K$ 
an inertia subgroup at~$\lambda$.

Before we prove the theorem, let us highlight the core idea.
We will work by contradiction, so
suppose there is an elliptic curve 
$E/K$ that is a counterexample to the theorem.
We show the existence of a quadratic twist
$F/K$ such that the $\ell$-adic Galois representation 
$\rho_{F,\ell} : G_K \rightarrow \GL_2(\Z_\ell)$
attached to $F$ satisfies 
\[
 \rho_{F,\ell}(G_K)=\rho_{F,\ell}(I_\lambda);
\]
that is, the global image of $\rho_{F,\ell}$ 
is equal to the image of a local Galois group. 
In general, a global image is huge whilst a local image is much smaller, so this situation is prone to a contradiction, which 
we will show to happen in our setting.

\smallskip

We now prove the theorem. Suppose there is an elliptic 
curve $E/K$ having a $K$-rational $\ell$-isogeny, good reduction
away from $\lambda$, and potentially multiplicative
reduction at~$\lambda$.

\smallskip

\textbf{Claim}:
there is a quadratic twist $F/K$
of $E$ of conductor $\lambda$
such that  
$K(F[\ell^n])/K$
is an $\ell$-extension for all $n \ge 1$.

\smallskip

We first show how this claim implies the theorem.
Indeed, by the 
criterion of 
N\'eron--Ogg--Shafarevich \cite[Proposition IV.10.3]{SilvermanII},
the extension $K(F[\ell^n])/K$ is unramified
away from $\lambda$ and the infinite places.
By assumption (iii), we have $\ell \nmid h_K^+$, hence
part~(b) of Theorem~\ref{thm:Iwasawa} to deduce
that $\lambda$ is totally ramified in $K(F[\ell^n])/K$.
Let $n \ge 1$ and consider
\[
\overline{\rho}_{F,\ell^n} : G_K \rightarrow \GL(F[\ell^n]) \cong \GL(\Z/\ell^n\Z)
\]
the mod $\ell^n$ representation of $F$.
The Galois group of the extension $K(F[\ell^n])/K$
is
$\overline{\rho}_{F,\ell^n}(G_K)$
and its inertia subgroup at $\lambda$
is $\overline{\rho}_{F,\ell^n}(I_\lambda)$.
As the extension is totally ramified, we have
\[\overline{\rho}_{F,\ell^n}(G_K)=
\overline{\rho}_{F,\ell^n}(I_\lambda).\]
However, the latter group is reducible 
by the theory of the Tate curve (see \cite[Exercise V.5.13]{SilvermanII}). 
As $F$ does not have complex multiplication (it has a multiplicative prime)
this contradicts Serre's open image theorem \cite[Chapter IV]{SerreBook}
for sufficiently large $n$. (Note that taking the limit 
on~$n$ leads to the equality of $\ell$-adic representations 
as in the discussion above.)

\smallskip

It remains to establish our claim.
The curve $E$ has potentially multiplicative
reduction at $\lambda$, and the theory of the Tate curve
(c.f.\ \cite[Exercises V.5.11 and V.5.13]{SilvermanII}) gives a precise 
description of the restriction of the representation $\overline{\rho}_{E,\ell}$
to 
the decomposition group $D_\lambda$:
\begin{equation}\label{eqn:inertiaTate}
\overline{\rho}_{E,\ell} \vert_{D_\lambda} \sim
\begin{pmatrix}
\eta \cdot \chi_{\ell} & * \\
0 & \eta 
\end{pmatrix}
\end{equation}
where $\chi_{\ell}$ is the modulo $\ell$ cyclotomic character,
and $\eta$ is a character of $D_\lambda$ which is trivial or quadratic.
Moreover, the (local) twist $E \otimes \eta$ is an
elliptic curve defined over $K_\lambda$ having split multiplicative
reduction at $\lambda$. 

As $E$ has a $K$-rational $\ell$-isogeny, the mod $\ell$ 
representation is reducible:
\[
\overline{\rho}_{E,\ell} \sim
\begin{pmatrix}
\phi & * \\
0 & \psi
\end{pmatrix}
\]
where $\phi$, $\psi$ are characters $G_K \rightarrow \F_\ell^*$.
It follows from the criterion 
of N\'eron--Ogg--Shafarevich
that $\phi$ and $\psi$
are unramified except possibly at $\lambda$ and the infinite places.
By assumption (i), the mod $\ell$ cyclotomic character is 
trivial on $G_K$. From \eqref{eqn:inertiaTate}
we have $\phi \vert_{I_\lambda}=\psi \vert_{I_\lambda}=\eta \vert_{I_\lambda}$ is of order dividing~2.
Thus $\phi/\psi$ and $\phi^2$ are characters of~$G_K$
of order dividing $\ell-1$ that are unramified
away from the infinite places. By assumption (iii),
the narrow class number $h_K^+$ is coprime to $\ell-1$,
thus $\phi/\psi=\phi^2=1$.
Hence $\phi=\psi$ is a quadratic character of $G_K$.
We let $F$ be the (global) quadratic twist
$E \otimes \phi$. Now $F/K$ has conductor $\lambda$,
and 
\[
\overline{\rho}_{F,\ell} \sim
\begin{pmatrix}
1 & * \\
0 & 1 
\end{pmatrix}.
\]
Thus $\#\overline{\rho}_{F,\ell}(G_K)=1$ or $\ell$, 
hence $\overline{\rho}_{F,\ell}(G_K)$ is an $\ell$-group.
To complete the proof of our claim
we need to show that $\overline{\rho}_{F,\ell^n}(G_K)$
is an $\ell$-group for all $n$.
Consider the commutative diagram
\[
\xymatrixcolsep{5pc}\xymatrix{
G_K \ar@{->}[rd]_{\overline{\rho}_{F,\ell}} \ar@{->}[r]^{\overline{\rho}_{F,\ell^n}}
&\GL_2(\Z/\ell^n\Z) \ar@{->}[d]^\pi\\
& \GL_2(\F_\ell)} 
\]
where $\pi$ is the reduction modulo $\ell$ map. 
From this we deduce the exact sequence
\[
1 \rightarrow \overline{\rho}_{F,\ell^n} (G_K) \cap \ker(\pi)
\rightarrow \overline{\rho}_{F,\ell^n} (G_K) \rightarrow \overline{\rho}_{F,\ell}(G_K)
\rightarrow 1.
\]
We know already that $\overline{\rho}_{F,\ell}(G_K)$ is an $\ell$-group.
Thus
it is sufficient to show that $\ker(\pi)$ is an $\ell$-group.
However,  
\[
\ker(\pi)=\left\{\begin{pmatrix}
a & b \\ 
c & d
\end{pmatrix} \; : \; a,~b,~c,~d \in \Z/\ell^n \Z, \quad
a \equiv d \equiv 1 , \; b \equiv c \equiv 0 \pmod{\ell}
\right\}.
\]
We see that $\# \ker(\pi)=\ell^{4n-4}$, proving our claim,
and completing the proof. \qed

\section{Criteria for Asymptotic FLT}\label{sec:criteria}
Freitas and Siksek \cite{FS} give a criterion
for the asymptotic FLT over a totally real field $K$
in terms of solutions to a certain $S$-unit equation.
\c{S}eng\"{u}n and Siksek \cite{Haluk} give similar criterion
for general number fields, assuming standard conjectures
that we state below. 
In this section we
state these results, and then show that they
can be sharpened for
number fields where $2$ is totally ramified.
\subsection*{The $S$-unit equation}
Let $K$ be a number field.
Let
\begin{equation}\label{eqn:ST} 
\begin{gathered}
S=\{ \mP \; :\; \text{$\mP$ is a prime of $K$ above $2$}\}, \\
T=\{ \mP \in S \; : \; 
f(\mP/2)=1\},
\qquad 
U=\{ \mP \in S \; : \; 
3 \nmid \ord_\mP(2) \}.
\end{gathered}
\end{equation}
Here $f(\mP/2)$ denotes the residual degree of $\mP$.
The ring of $S$-integers is 
\[
\OO_S=\{ \alpha \in K \; : \; \text{$\ord_\fq(\alpha) \ge 0$
for all primes $\fq \notin S$}\}.
\]
We let $\OO_S^\times$ be the unit group of $\OO_S$; explicitly
\[
\OO_S^\times=\{\alpha \in K \; : \; \text{$\ord_\fq(\alpha) = 0$
for all primes $\fq \notin S$}\} .
\]
Consider the $S$-unit equation
\begin{equation}\label{eqn:sunit}
\lambda+\mu=1, \qquad \lambda,~\mu \in \OO_S^\times,
\end{equation}
and let $\sS_3$ be the subgroup of $\PGL_2(K)$ given by
\[
\sS_3=\{\, z, \, 1/z,\, 1-z, \, 1/(1-z),\, z/(z-1),\, (z-1)/z\}.
\]
As explained in \cite[Section 6]{FS}, there is an action of 
$\sS_3$ on the solutions to \eqref{eqn:sunit} given
by $\sigma(\lambda,\mu)=(\sigma(\lambda),1-\sigma(\lambda))$
for $\sigma \in \sS_3$.
The solutions $(2,-1)$, $(-1,2)$, $(1/2,1/2)$ form
an orbit under the action of $\sS_3$ that we call the
\textit{irrelevant orbit}. Other solutions are called \textit{relevant}.

\subsection*{A criterion of asymptotic FLT over totally real fields}
The following is (a special case of) Theorem~3 of \cite{FS}.
\begin{thm}\label{thm:FS}
Let $K$ be a totally real field
and let $S$, $T$ and $U$ be as  \eqref{eqn:ST}.  Suppose 
that either $T \ne \emptyset$ or $[K:\Q]$ is odd.
Suppose that for every solution $(\lambda,\mu)$ to the $S$-unit 
equation~\eqref{eqn:sunit}
there is
\begin{enumerate}
\item[(A)] either 
some $\mP \in T$ that satisfies
$\max\{ \lvert \ord_{\mP} (\lambda) \rvert, \lvert \ord_{\mP}(\mu) \rvert \} 
\le 4 \ord_{\mP}(2)$.
\item[(B)] or some $\mP \in U$ that satisfies both
$\max\{ \lvert \ord_{\mP} (\lambda) \rvert, \lvert \ord_{\mP}(\mu) \rvert \} 
\le 4 \ord_{\mP}(2)$, and
$\ord_{\mP}(\lambda \mu) \equiv \ord_{\mP}(2) \pmod{3}$.
\end{enumerate}
Then the asymptotic FLT holds over~$K$.
\end{thm}

\subsection*{A criterion of asymptotic FLT over general fields}
For now $K$ will be a general number field---by general
we simply mean that we do not require it to be totally real. We will soon state
a theorem of \c{S}eng\"{u}n and Siksek \cite{Haluk}
which gives a criterion for asymptotic FLT over general number fields.
This criterion is dependent on two standard conjectures which we now
state in precise form. However we do not apply the conjectures
directly, and so we omit any further discussion of them; 
instead we recommend
the exposition of these conjectures
in \cite[Sections 2--4]{Haluk} as well the references
cited therein.
The first conjecture is a special case 
of Serre's modularity conjecture over $K$.
\begin{conj} \label{conj:Serre}
Let $\overline{\rho} : G_K \rightarrow \GL_2(\overline{\F}_p)$ be an odd,
irreducible, continuous representation with Serre conductor $\cN$ (prime-to-$p$
part of its Artin conductor) and such that
$\det(\overline{\rho})=\chi_p$ is the mod~$p$
cyclotomic character.
Assume that $p$ is unramified in~$K$ and that
$\overline{\rho}\vert_{G_{K_\fp}}$ arises from a finite-flat group scheme over
$\OO_\fp$ for every prime~$\fp \mid p$. Then there is a (weight $2$) mod $p$
eigenform $\theta$ over $K$ of level
$\cN$ such that for all primes $\fq$ coprime to $p\cN$, we have
\[
\Tr(\overline{\rho}
({\Frob}_\fq)) = \theta(T_\fq),
\]
where $T_\fq$ denotes the Hecke operator at $\fq$.
\end{conj}
\noindent We point out that the condition $\det(\overline{\rho})=\chi_p$ is
inadvertently omitted in \cite{Haluk}.

\smallskip

The second conjecture is in essence a generalization of the 
Eichler--Shimura theorem for rational weight $2$ eigenforms.
\begin{conj}\label{conj:ES}
  Let $\ff$ be a (weight $2$) complex eigenform over $K$ of level $\cN$ 
that is non-trivial and new. 
If $K$ has some real place, then there exists an elliptic curve $E_\ff/K$, of conductor
  $\cN$, such that
 \begin{equation}\label{eq:pd}
    \# E_\ff(\OO_K/\fq) = 1 + {\bf N}\fq - \ff(T_\fq)
\quad \text{for all $\fq \nmid \cN$}.
\end{equation}
If $K$ is totally complex, 
then there exists either an elliptic curve $E_\ff$ of conductor $\cN$ 
satisfying \eqref{eq:pd} 
or a fake elliptic curve $A_\ff/K$, of conductor
$\cN^2$, such that
\[
    \# A_\ff(\OO/\fq) = \left (1 + {\bf N}\fq - \ff(T_\fq) \right )^2
\quad \text{for all $\fq \nmid \cN$}.
\]
\end{conj}
The following is Theorem 1.1 of \cite{Haluk}.
\begin{thm}\label{thm:Haluk}
Let $K$ be a number field over which Conjectures~\ref{conj:Serre}
and~\ref{conj:ES} hold.
Let $S$, $T$ be as in \eqref{eqn:ST}.
Suppose $T \ne \emptyset$.
Suppose that for every solution $(\lambda,\mu)$ to the $S$-unit 
equation~\eqref{eqn:sunit}
there is
some $\mP \in T$ that satisfies
$\max\{ \lvert \ord_{\mP} (\lambda) \rvert, \lvert \ord_{\mP}(\mu) \rvert \} 
\le 4 \ord_{\mP}(2)$.
Then the asymptotic FLT holds over~$K$.
\end{thm}

\subsection*{A simplification when $\#S=1$}
The following proposition is a simplification and strengthening of
the ideas in \cite[Section 6]{FS}, under the additional
hypothesis that the set~$S$ has precisely one element.
It allows us to simplify the condition on the $S$-unit equation
in Theorems~\ref{thm:FS} and~\ref{thm:Haluk}. 
The proof makes use of ideas found in 
the proof of \cite[Lemme 1]{Kraus}.
\begin{prop}\label{prop:twist}
Let $K$ be a number field with precisely one prime $\mP$
above $2$, and let $S=\{\mP\}$.
The following are equivalent.
\begin{enumerate}
\item[(i)] There is an elliptic curve $E/K$ with full $2$-torsion
and conductor $\mP$.
\item[(ii)] There is an elliptic curve $E/K$ with full $2$-torsion,
potentially good reduction away from $\mP$, and potentially
good multiplicative reduction at $\mP$.
\item[(iii)] There is a solution $(\lambda,\mu)$
to the $S$-unit equation \eqref{eqn:sunit}
with 
$\ord_\mP(\mu)=0$ and $\ord_\mP(\lambda) > 4 \ord_\mP(2)$.
\item[(iv)] There is a solution $(\lambda,\mu)$
to the $S$-unit equation \eqref{eqn:sunit}
with 
\[
\max\{\lvert \ord_\mP(\lambda)\rvert, \lvert \ord_\mP(\mu) \rvert\}
> 4 \ord_\mP(2).
\]
\end{enumerate}
\end{prop}
\begin{proof}
Clearly (iii) implies
(iv). If we have a solution 
as in (iv), we simply observe that we can produce a solution
satisfying (iii) by applying a suitable element of $\sS_3$.
This proves the equivalence of (iii) and (iv). 
Clearly (i) implies (ii). 
To complete the proof it is 
enough to show that (ii) implies (iii) and that (iii) implies (i).

We suppose (ii).
The elliptic curve $E$ has the form
\[
E \; : \; Y^2=X(X-a)(X+b)
\]
with $a$, $b \in K$ and $ab(a+b) \ne 0$. Let $c=-a-b$.
Then $a+b+c=0$. Applying a permutation
to $a$, $b$, $c$ allows us to
suppose that $\ord_\mP(b) \ge \ord_\mP(c) \ge \ord_\mP(a)$.
If this permutation is cyclic then the resulting
elliptic curve is isomorphic to our original
model, and if non-cyclic then it is a quadratic
twist by $-1$.

Let $\lambda=-b/a$, $\mu=-c/a$. Then,
\begin{equation}\label{eqn:presunit}
\lambda+\mu=1, \qquad \ord_\mP(\lambda) \ge \ord_\mP(\mu) \ge 0.
\end{equation}
 The quadratic twist of $E$ by $-a$ is 
\begin{equation}\label{eqn:Eprime}
E^\prime \; : \; Y^2=X(X+1)(X+\lambda).
\end{equation}
and also has potentially multiplicative reduction at $\mP$ and 
potentially good reduction away from $\mP$.
In the usual notation, the invariants of $E^\prime$ are
\[
\begin{gathered}
c_4=16(\lambda^2-\lambda+1), \qquad c_6=-64(1-\lambda/2)(1-2\lambda)(1+\lambda)
\\
\Delta=16\lambda^2(\lambda-1)^2, \qquad j=\frac{2^8 (\lambda^2-\lambda+1)^3}{\lambda^2(\lambda-1)^2}.
\end{gathered}
\]
If $\fq$ is a prime $\nmid 2$ then $\ord_\fq(j) \ge 0$, and we easily check from the above
formulae that this forces $\ord_\fq(\lambda) = \ord_\fq(\mu)=
\ord_\fq(\Delta)=0$. 
If $\ord_\mP(\lambda)=0$ then, by \eqref{eqn:presunit},
 $\ord_\mP(1-\lambda)=0$
and so $\ord_\mP(j)>0$ giving a contradiction. 
Thus $\ord_\mP(\lambda)=t$ with $t>0$, and as $\mu=1-\lambda$,
we have $\ord_\mP(\mu)=0$. 
Since 
$\ord_\mP(j)=8\cdot \ord_\mP(2)-2t<0$ we have $\ord_\mP(\lambda)=t>4\ord_\mP(2)$. 
This proves that (ii) implies (iii).

We now suppose $(\lambda,\mu)$ are as in (iii).
We let $E^\prime$ be the elliptic curve given by \eqref{eqn:Eprime}.
This model is integral and has good reduction at all $\fq \nmid 2$.
As $\ord_\mP(\lambda) > 4 \ord_\mP(2)$,
 Hensel's Lemma shows that the expressions 
$\lambda^2-\lambda+1$, $1-\lambda/2$, $1-2 \lambda$
and $1+\lambda$ are all $\mP$-adic squares. Thus $-c_4/c_6$ is an $\mP$-adic square. By \cite[Theorem V.5.3]{SilvermanII}
the elliptic curve $E^\prime$ has split multiplicative reduction at $\mP$. 
This shows that (iii) implies (i).
\end{proof}

\begin{thm}\label{thm:llsimple}
Let $K$ be a number field. 
If $K$ has complex embeddings, assume Conjectures~\ref{conj:Serre}
and~\ref{conj:ES}.
Suppose $2$ is totally ramified in $K$
and let $\mP$ be the prime above $2$. Suppose there is no
elliptic curve $E/K$ with full $2$-torsion and conductor~$\mP$.
Then the asymptotic FLT holds over~$K$.
\end{thm}
\begin{proof}
Here $S=T=\{\mP\}$. By Proposition~\ref{prop:twist}, 
as there is no elliptic curve 
with full $2$-torsion and conductor $\mP$,
every solution to \eqref{eqn:sunit} satisfies
$\max\{\lvert \ord_\mP(\lambda) \rvert,\lvert \ord_\mP(\mu) \rvert\} \le 4 \ord_\mP(2)$. 
The theorem follows from Theorem~\ref{thm:FS} if $K$ is totally real,
and from Theorem~\ref{thm:Haluk} otherwise.
\end{proof}

\section{Proof of Theorem~\ref{thm:main}}
Let $K$ be a number field.
The narrow Hilbert $p$-class field of $K$ is the maximal abelian $p$-extension
of $K$ unramified away from the finite places. This is a subextension
of the narrow Hilbert class field of degree $p^{\ord_p(h_K^+)}$. 
We thank D.\ Bernardi for suggesting the proof of the following lemma.
\begin{lem}\label{lem:Bernardi}
Let $K$ be a number field and $L$ the narrow Hilbert $p$-class field of $K$.
Suppose $\Gal(L/K)$ is cyclic. Then 
$p \nmid h_L^+$.
\end{lem}
\begin{proof}
Let $M$ be the narrow Hilbert $p$-class field of $L$. 
To establish the lemma we merely have to show that $M=L$. 
First we claim that $M/K$ is Galois. To see this suppose
that $\sigma \in \Gal(\overline{K}/K)$ where we identify
$L$, $M$ as subfields of $\overline{K}$. 
Now $L/K$ is Galois and so $L^\sigma=L$. Thus $M$
and $M^\sigma$ are both the unique maximal abelian $p$-extension of $L=L^\sigma$
unramified at all finite places, and so must be equal. It follows that $M/K$ is Galois. We remark that $M/K$ is 
a $p$-extension (Lemma~\ref{lem:galClosure}) unramified at all finite places.

Let $G=\Gal(M/K)$ and let $G^\prime$ its derived subgroup. 
As
$\Gal(L/K) \cong G/\Gal(M/L)$ is abelian, the subgroup
$\Gal(M/L)$ contains $G^\prime$. 
Thus the fixed field $M^{G^\prime}$ satisfies 
$L \subset M^{G^\prime} \subset M$. 
However $\Gal(M^{G^\prime}/K) \cong G/G^\prime$
which is abelian; thus $M^{G^\prime}$ is an abelian $p$-extension
of $K$ unramified at the finite places and so is contained in $L$.
It follows that $M^{G^\prime}=L$. Hence $\Gal(M/L) \cong G^\prime$
and $G/G^\prime \cong \Gal(L/K)$ which is cyclic.
From Corollary~\ref{cor:abIsCyclic}, 
the derived group $G^\prime \cong \Gal(M/L)$ is trivial,
and so $M=L$ as required. 
\end{proof}
\begin{prop}\label{prop:equiv}
Let $K$ is a number field with exactly one
prime $\mP$ above $2$. 
Let $L$ be the narrow Hilbert $2$-class field of $K$.
The following are equivalent.
\begin{enumerate}
\item[(i)] There exists a finite extension $M/K$
with unique prime above $2$ and odd narrow class number.
\item[(ii)] $L$ has unique prime above $2$.
\item[(iii)] $L$ has a unique prime above $2$
and odd narrow class number.
\item[(iv)] The order of $[\mP]$ in $\Clp(K)$
is divisible by $h_{K,2}^+$.
\end{enumerate}
\end{prop}
\begin{proof}
Suppose (i) is satisfied. Then $LM$ is a $2$-extension of $M$
unramified away from the infinite places. Thus $LM=M$,
and so $L \subseteq M$. In particular, $L$ has a unique prime
above $2$.  Thus (ii) is satisfied.

Now suppose (ii) is satisfied. As $\mP$ is unramified in $L$,
and there is a unique prime above $\mP$, we see that $\mP$
is inert in $L$. Thus $\Gal(L/K)$ is equal to the decomposition
group at $\mP$ which is cyclic. By Lemma~\ref{lem:Bernardi}
we know that $L$ has odd narrow class number and so (iii) is
satisfied. Clearly, (iii) implies both (i) and (ii) so 
we have now proved the equivalence of (i), (ii), (iii).

We complete the proof by showing that (ii) and (iv) are equivalent.
Let $N$ be the narrow class field of $K$. We consider the Artin
map for the extension $N/K$. This is an
isomorphism $\theta : \Cl^+(K) \rightarrow \Gal(N/K)$, 
sending $[\mP]$ to the Frobenius element $\Frob_{N/K,\mP}$.
We compose this with the natural restriction
map $\phi: \Gal(N/K) \twoheadrightarrow 
\Gal(L/K)$. Let $H$ be the $2$-Sylow subgroup of $\Gal(N/K)$,
and $H^\prime$ be the maximal subgroup of odd order.
As $\Gal(N/K)$ is abelian, $\Gal(N/K)=H\oplus H^\prime$.
Recall that $\#\Gal(L/K)=h_{K,2}^+=\#H$.
Thus the restriction
$\phi\vert_{H} : H \rightarrow \Gal(L/K)$
is an isomorphism.
The map $\phi$ sends $\Frob_{N/K,\mP}$ to $\Frob_{L/K,\mP}$.
Note that $L$ has a unique prime above $2$ if and only if $\Frob_{L/K,\mP}$ has
order $\#\Gal(L/K)=h_{K,2}^+$.  This is equivalent to $\Frob_{N/K,\mP}$
having order divisible by $h_{K,2}^+$, which is equivalent to $[\mP]$
having order divisible by $h_{K,2}^+$.
\end{proof}

\begin{cor}\label{cor:KrausCrit}
Let $K$ is a number field with exactly one
prime $\mP$ above $2$. Suppose that the
order $[\mP]$ in $\Clp(K)$
is divisible by $h_{K,2}^+$. Then
there is no elliptic curve $E/K$
with a $K$-point of order $2$, good reduction
away from $\mP$, and potentially multiplicative
reduction at $\mP$.
\end{cor}
\begin{proof} 
Let $L$ be the narrow $2$-class field of $K$.
By Proposition~\ref{prop:equiv} the number field $L$
has one prime $\lambda$ above $2$ and odd narrow class number.
Suppose there is an elliptic curve $E/K$
as in the statement. Recall that good reduction and potentially 
multiplicative reduction are preserved by base change.
Therefore, the curve 
$E/L$ has an $L$-point of order $2$, good reduction
away from $\lambda$, and potentially multiplicative
reduction at~$\lambda$.
Now Theorem~\ref{thm:Krausgen}
applied over~$L$  with $\ell=2$ gives a contradiction.
\end{proof}

\subsection{Proof of Theorem~\ref{thm:main}} 
Theorem~\ref{thm:main} follows immediately from
Theorem~\ref{thm:llsimple} 
together with Corollary~\ref{cor:KrausCrit}.

\section{Proof of Theorem~\ref{thm:main2}}
In this section $K$, $\mP$ satisfy assumptions (a), (b), (c)
of Theorem~\ref{thm:main2}. Moreover, $K^\prime/K$ is a $2$-extension
that is unramified away from $\mP$ (and in particular,
unramified at the infinite places).
By the first part of Theorem~\ref{thm:Iwasawa},
the prime $\mP$ is totally ramified in $K^\prime$, and thus $2$ is totally
ramified in $K^\prime$. We stress that it is here
that we have made use of the fact that $h_K$ is odd (assumption (c)). 
We denote by $\mP^\prime$ the unique prime of $K^\prime$ above~$2$.

Next we let $L$ be the narrow $2$-class field of $K$. As assumptions (a),
(b) are satisfied, Proposition~\ref{prop:equiv} tells that 
$h_{L}^+$ is odd, and $L$ has unique prime $\fq$ (say) above $2$.
Let $L^\prime=L K^\prime$.
As $K^\prime/K$ is a $2$-extension unramified away from $\mP$,
the extension $L^\prime /L$ is a $2$-extension unramified
away from $\fq$. 
Now we apply the second part of Theorem~\ref{thm:Iwasawa}
to the extension $L^\prime /L$ to deduce that $h_{L^\prime}^+$
is odd, and that there is a unique prime $\fq^\prime$ say above $\fq$
(and so $\fq^\prime$ is the unique prime of $L^\prime$ above $2$).
In particular, by Theorem~\ref{thm:Krausgen}, there is no 
elliptic curve $E/L^\prime$ with an $L^\prime$-rational
point of order $2$ and conductor $\fq^\prime$. 
As $L^\prime \supseteq K^\prime$,
it follows (as in the proof of Corollary~\ref{cor:KrausCrit}) that there is no elliptic curve $E/K^\prime$
with $K^\prime$-rational point of order $2$
and conductor~$\mP^\prime$.

The proof is completed by applying Theorem~\ref{thm:llsimple} with $K^\prime$
and $\mP^\prime$ instead of $K$ and $\mP$. 
  \qed

\subsection*{Example}
We give an example to illustrate the importance of assumption (c) for proof 
of Theorem~\ref{thm:main2}.
Let $K=\Q(\alpha)$ be the totally real quartic
field with $\alpha$ satisfying $\alpha^4 - 10\alpha^2 - 8\alpha + 7 = 0$.
Then $2$ totally ramifies in $K$ and we denote the prime
above $2$ by $\mP$. We checked using the computer algebra
package \texttt{Magma} that $\Cl^+(K)$ is cyclic
of order $4$ generated by $[\mP]$. Thus asymptotic FLT holds
for $K$ by Theorem~\ref{thm:main}. 

However $h_K=2$, and so condition (c) is not satisfied.
The Hilbert class field of $K$ is in fact the quadratic extension
$K(\sqrt{2})=K \cdot \Q(\zeta_8)^+$. 
Of course $\mP$ is 
not ramified 
in this extension, and so by Proposition~\ref{prop:equiv}
 it is inert. 
It follows that the unique prime $\mP_r$ above $2$ 
in $K_r=K \cdot \Q(\zeta_{2^r})^+$
has residue field $\F_4$ for all $r$.
We are therefore
unable to apply Theorem~\ref{thm:llsimple}
(with $K_r$ in place of $K$). By 
\cite[Theorem~3]{FS} if there is
an non-trivial solution to the Fermat equation over $K_r$
with $p$ suitably large, and we allow ourselves to
assume a suitable Eichler--Shimura conjecture,
 then there is an elliptic curve $F/K_r$
with full $2$-torsion and good reduction away from $\mP_r$. However, this
result does not specify potentially multiplicative reduction
at $\mP_r$ when the residue field is not $\F_2$.

\section{Proof of Corollary~\ref{cor:correct}}\label{sec:quadratic}
In this section $K$ is a quadratic field (real or imaginary).
We want to understand when $K$ satisfies conditions (a), (b), (c)
of Theorem~\ref{thm:main2}.
We start by recalling some basic facts
from genus theory, following \cite[Section 1.3.1]{Lemmermeyer}.
Let $D$ be the discriminant of $K$. 
The \textbf{prime discriminants} are the integers $-4$, $8$, $-8$,
and $(-1)^{(\ell-1)/2} \ell$, where $\ell$ is an odd prime.
The discriminant $D$ can be written as a product $D=d_1 d_2 \cdots d_t$
where the $d_j$ are prime discriminants, and this factorization is
unique up to reordering. The $2$-ranks of the class group
and the narrow class group of $K$ have convenient expressions
in terms of this factorization:
\begin{equation}\label{eqn:genus}
\dim_{\F_2} \Cl^+(K)[2]=t-1,
\quad
\dim_{\F_2} \Cl(K)[2]=
\begin{cases}
t-1 & \text{if $D<0$}\\
t-1 & \text{if all $d_j>0$,}\\
t-2 & \text{if $D>0$ and some $d_j<0$.}
\end{cases}
\end{equation}
For the remainder of this section we assume
that $2$ ramifies in $K$, i.e. condition~(a). 
After reordering the $d_j$ if necessary,
we have $d_1=-4$, $8$ or $-8$. The following is immediate
from \eqref{eqn:genus}.
\begin{lem}\label{lem:oddnarrow}
$h_K^+$ is odd if and only if $K=\Q(i)$, $K=\Q(\sqrt{2})$ or $K=\Q(\sqrt{-2})$.
\end{lem}
\begin{lem}\label{lem:oddeven}
$h_K$ is odd and $h_K^+$ is even if and only if $K=\Q(\sqrt{\ell})$
or $K=\Q(\sqrt{2\ell})$ where $\ell \equiv 3 \pmod{4}$ is a prime.
In this case the $2$-Sylow subgroup of $\Cl^+(K)$ has order $2$.
\end{lem} 
\begin{proof} 
From ~\eqref{eqn:genus} we see that $h_K$ is odd and $h_K^+$ is even if and only if
$t=2$, $d_1<0$, and $d_2<0$. In particular, $d_1 = -4, -8$. 
Recall that the discriminant $D$ of a quadratic field $\Q(\sqrt{d})$ ($d$ square-free) is $d$ if $d \equiv 1 \pmod{4}$ and $4d$ otherwise. 
From this it is easy to deduce
the first part of the lemma. 

It then follows in this case that 
$\Cl^+(K)[2] \cong \Z/2\Z$. Write $\Cl^+(K)[2^\infty]$
for the $2$-Sylow subgroup of $\Cl^+(K)$. 
For the second part of the lemma
we want to show
that $\Cl^+(K)[2^\infty]=\Cl^+(K)[2]$. This can easily
be deduced from a standard exact sequence \cite[Proposition 3.2.3]{Cohen}
relating the class group and narrow class group. Alternatively,
let $\ga$ be a non-zero ideal of $\OO_K$
representing a class $[\ga]$ in $\Cl^+(K)[2^\infty]$;
we would like to show that $[\ga] \in \Cl^+(K)[2]$. As $h_K$ is odd,
$\ga^r=\alpha \cdot \OO_K$ for some non-zero  $\alpha \in \OO_K$
with $r$ odd. Then $\ga^{2r} = \alpha^2 \OO_K$ and $\alpha^2$
is totally positive, and so $[\ga]^{2r}=1$ in $\Cl^+(K)$. 
As $r$ is odd, $[\ga] \in \Cl^+(K)[2]$.
\end{proof}

\subsection*{Application of Theorem~\ref{thm:main2} to Quadratic Fields}
We shall prove Corollary~\ref{cor:correct}. The second part
of the corollary follows from the first part and
Theorem~\ref{thm:main2} as the fields $K_r$ are totally 
real $2$-extensions of $K$.

We now prove the first part.
For this we would like to know all real quadratic
fields satisfying conditions (a), (b), (c) in the
statement of Theorem~\ref{thm:main2}. 
By Lemmas~\ref{lem:oddnarrow}~and~\ref{lem:oddeven} we see these are the fields
$\Q(\sqrt{2})$ and
$\Q(\sqrt{d})$ where $d=\ell$ or $d=2\ell$ and $\ell \equiv 3 \pmod{4}$ is a prime. Moreover, from the second part of Lemma~\ref{lem:oddeven}, we see that to prove the corollary it is enough to show that $[\mP]$
generates $\Cl^+(K)[2^\infty]=\Cl^+(K)[2]\cong \Z/2\Z$
if and only if $\ell \equiv 3 \pmod{8}$. 
As $2$ ramifies in $K$, we have $\mP^2=2\OO_K$.
Thus $[\mP]$ has order $1$ or $2$ in $\Cl^+(K)$
and it is precisely in the latter case that $[\mP]$
generates $\Cl^+(K)[2]$.
However the class number is odd, so $\mP$ must
be principal. We may write $\mP=\beta \OO_K$, with $\beta=a+b\sqrt{d}>0$
where $a$, $b \in \Z$. 
Since $\ell \equiv 3 \pmod{4}$, quadratic reciprocity implies
the norm of any unit is $1$, and so every unit is totally positive
or totally negative. Our positive generator $\beta$ for $\mP$
is therefore unique up to multiplication
by a totally positive unit. Thus $[\mP]$ has order $2$ in $\Cl^+(K)$
if and only if $\overline{\beta}<0$, where $\overline{\beta}=a-b\sqrt{d}$.
However,
\[
\beta \cdot \overline{\beta}=a^2-d b^2=2 \cdot \eta, \qquad \eta=\pm 1.
\]
Thus $[\mP]$ has order $2$ in $\Cl^+(K)$ if and only if $\eta=-1$.
Now $\ell \mid d$, and so $2 \cdot \eta$ is a quadratic residue
modulo $\ell$. We deduce that 
$\eta=-1$ if and only if $\ell \equiv 3 \pmod{8}$ 
completing the proof.

\section{Proof of Theorem~\ref{thm:cubic}}\label{sec:cubic}
In this section $K$ is a totally real cubic field with
$2$ either totally ramified or inert, 
$3$ ramified, and $\Delta_K$ (the discriminant)
a non-square. We would like to show that $K$ satisfies
asymptotic FLT. For this we will apply Theorem~\ref{thm:FS}.
Write $\mP$ for the unique prime above $2$.
In the notation of that theorem (c.f. \eqref{eqn:ST}),
we have 
\[
\begin{cases}
S=T=\{\mP\},~U=\emptyset & \text{\textbf{case (I)}: if $2$ is totally ramified in $K$},\\
S=U=\{\mP\},~T=\emptyset & \text{\textbf{case (II)}: if $2$ inert in  $K$}.
\end{cases}
\]
We shall show
that conditions (A), (B) of Theorem~\ref{thm:FS}
respectively hold for all solutions
$(\lambda,\mu)$
to the $S$-unit equation \eqref{eqn:sunit} according
to whether we are in case (I) or (II).
Note that the action of
$\sS_3$ on 
$(\lambda,\mu)$ preserves the value 
$\max\{\lvert \ord_{\mP}(\lambda)\rvert, \lvert \ord_{\mP}(\mu)\rvert\}$,
 and also
the residue class of $\ord_\mP(\lambda\mu)$ modulo $3$. Thus
need only show that conditions (A), (B) hold for a representative
of each $\sS_3$-orbit. Now for a solution $(\lambda,\mu)$
we may apply a suitable element of $\sS_3$ so that
\[
\ord_\mP(\lambda)=m \ge 0, \qquad \ord_\mP(\mu)=0.
\]
Write
\[
\Norm(\lambda)=\eta_1 \cdot 2^n, \qquad \Norm(\mu)=\eta_2,
\qquad 
\eta_1=\pm 1,~\eta_2=\pm 1, \quad
n=\begin{cases}
m & \text{case (I)}\\
3m & \text{case (II)}.
\end{cases}
\]
If $\mu \in \Q$, then $\mu$ is a unit in $\Z$ and so
$\mu=\pm 1$ which gives $(\lambda,\mu)=(2,-1)$ which
satisfies (A), (B) respectively for cases (I), (II).
We therefore suppose $\mu \notin \Q$ and so $K=\Q(\mu)$.
The minimal polynomial of $\mu$ has the form
\[
f(X)=X^3+aX^2+bX-\eta_2,
\]
for some $a$, $b \in \Z$. Write $\Delta_f$ for the discriminant of $f$.
Then
\[
\Delta_f=[\OO_K : \Z[\mu] ]^2 \cdot \Delta_K.
\]
Recall that $3 \mid \Delta_K$ and $\Delta_K$ is not a square.
Therefore $3 \mid \Delta_f$ and $\Delta_f$ is not a square.

Note that the minimal polynomial for $\lambda=1-\mu$ is $-f(1-X)$ as $\lambda \not\in \Q$;
therefore its constant coefficient must be $-\Norm(\lambda)$.
We deduce
\[
b=\eta_1 \cdot 2^n-1+\eta_2-a,
\qquad
f=X^3+aX^2+(\eta_1 \cdot 2^n-1+\eta_2-a) X-\eta_2.
\]
The discriminant $\Delta_f$ is now an expression
that depends on only $a$, $n$, $\eta_1$, $\eta_2$;
we denote this by $\Delta(a,n,\eta_1,\eta_2)$.
Let $a_0 \in \{0,1,2\}$, $n_0 \in \{0,1\}$
satisfy $a \equiv a_0 \pmod{3}$,
$n \equiv n_0 \pmod{2}$.
Note that $2^n \equiv 2^{n_0} \pmod{3}$.
Thus
\[
\Delta(a_0,n_0,\eta_1,\eta_2) \equiv \Delta(a,n,\eta_1,\eta_2)=\Delta_f
\equiv 0 \pmod{3}.
\]
We computed $\Delta(a_0,n_0,\eta_1,\eta_2)$
for the $24$ possible $(a_0,n_0,\eta_1,\eta_2)$
with $a_0 \in \{0,1,2\}$, $n_0\in \{0,1\}$ and $\eta_1$, $\eta_2 \in \{1,-1\}$.
We found $\Delta(a_0,n_0,\eta_1,\eta_2) \equiv 0 \pmod{3}$
for precisely the following two possibilities
\begin{equation}\label{eqn:possibility}
(a_0,n_0,\eta_1,\eta_2)=(0, 0, -1, -1) \quad \text{or} \quad (0, 1, 1, -1).
\end{equation}
In particular $\eta_2=-1$. 

\smallskip

\noindent \textbf{Case (I).} 
We will show condition (A) of Theorem~\ref{thm:FS} holds, 
that is $m \le 12$. In fact we prove the stronger $m \le 5$.
Thus we suppose $m \ge 6$. In particular $\mu=1-\lambda
\equiv 1 \pmod{4 \OO_K}$ and so $-1=\eta_2=\Norm(\mu) \equiv 1 \pmod{4}$
giving a contradiction.

\smallskip

\noindent \textbf{Case (II).} In this case $\mP=2\OO_K$.
Thus $\mu=1-\lambda \equiv 1 \pmod{2^m \OO_K}$, and so
$-1=\eta_2=\Norm(\mu) \equiv 1 \pmod{2^m}$. Thus $m=0$ or $1$.
If $m=1$ then (B) is satisfied. So suppose $m=0$, and so $n=0$,
and thus $n_0=0$. From \eqref{eqn:possibility}
we deduce
$\eta_1=-1$. Hence
\[
f=X^3+aX^2-(a+3)X+1.
\]
We find that
\[
\Delta_f=
    (a^2 + 3a + 9)^2.
\]
This contradicts the fact that $\Delta_f$ is not a square,
and completes the proof of Theorem~\ref{thm:cubic}.

\section{Proof of Theorem~\ref{thm:compCrit}}\label{sec:compCrit}
We merely have to prove (A). Parts (B), (C) follow
from (A) and Theorem~\ref{thm:llsimple}. 
Before proving (A) we will take a closer look at $U$. 
We suppose $U$, $V$ satisfy the hypotheses of the theorem:
the index $[\OO_K^\times : V]$ is finite and odd, and every element of 
$U$ is a square.
Let
\[
W:=\{ w \in \OO_K^\times \; : \; w \equiv 1 \pmod{16 \mP}\}.
\]
Then $U$ is contained in $W$ and we claim that the index $[W:U]$
is finite and odd. Indeed, $W$ is the kernel of the 
natural map $\OO_K^\times \rightarrow (\OO_K/16 \mP)^\times$,
and $U$ is the kernel of the restriction of this map to 
$V$. Consider the commutative diagram
\[
\xymatrix{
1 \ar[r]  & U \ar@{^{(}->}[d]
 \ar[r]  & V \ar[r] \ar@{^{(}->}[d] & (\OO_K/16\mP)^\times \ar@{=}[d] \ar[r] & 1 \\
1 \ar[r] & W  \ar[r] & \OO_K^\times \ar[r] & (\OO_K/16\mP)^\times \ar[r] & 1\\
} 
\]
The snake lemma immediately gives $W/U \cong \OO_K^\times/V$,
and thus the index $[W:U]$ is finite and odd. 

Next we show that every element of $W$ is a square. Let $(\OO_K^\times)^2$
be the subgroup of squares in $\OO_K^\times$. The assumption that 
every element in $U$ is a square is equivalent to saying that
$U$ is contained in $(\OO_K^\times)^2$. As the index 
$[\OO_K^\times : (\OO_K^\times)^2]$ is a power of $2$, and the
index $[W:U]$ is odd, we see that $W$ is also contained in $(\OO_K^\times)^2$.
We have now established our claim that every element of $W$ is a square.

We turn to the proof of  (A). Suppose (A) is false. 
By Proposition~\ref{prop:twist}
there is a solution $(\lambda,\mu)$ to the $S$-unit equation
satisfying
\begin{equation}\label{eqn:tlm}
\ord_\mP(\mu)=0, \qquad \ord_\mP(\lambda) >4 \ord_{\mP}(2).
\end{equation}
As there are only finitely many solutions to the $S$-unit equation,
we may suppose that $(\lambda,\mu)$ satisfies \eqref{eqn:tlm}
with the value of $\ord_{\mP}(\lambda)$ as large as possible.
Observe that $\mu \in \OO_K^\times$.
Moreover,
$\mu=1-\lambda$ and so $\mu \equiv 1 \pmod{16 \mP}$. Thus $\mu \in W$.
It follows that $\mu=\varepsilon^2$ for some $\varepsilon \in \OO_K^\times$.
We may therefore rewrite
\eqref{eqn:sunit} as $(1+\varepsilon)(1-\varepsilon)=\lambda$.
Hence
\begin{equation}\label{eqn:newrels}
1+\varepsilon=\lambda_1,\qquad 1-\varepsilon=\lambda_2, \qquad \lambda_1\lambda_2=\lambda.
\end{equation}
Here $\lambda_1$, $\lambda_2$ are in $\OO_K \cap \OO_S^\times$. Moreover,
by interchanging $-\varepsilon$ with $\varepsilon$ 
if necessary,
we may suppose that
\[
\ord_\mP(\lambda_1) \ge \ord_\mP(\lambda_2).
\]
However $\ord_\mP(\lambda_1)+\ord_\mP(\lambda_2)=\ord_\mP(\lambda)>4 \ord_\mP(2)$.
Thus $\ord_\mP(\lambda_1)> 2  \ord_\mP(2)$. Now from \eqref{eqn:newrels}
we have
\begin{equation}\label{eqn:prenew}
2=\lambda_1+\lambda_2, \qquad 2 \varepsilon=\lambda_1-\lambda_2.
\end{equation}
We immediately deduce that
\[
\ord_\mP(\lambda_2)=\ord_\mP(2), \qquad \ord_\mP(\lambda_1)=\ord_\mP(\lambda)-\ord_\mP(\lambda_2)=
\ord_\mP(\lambda)-\ord_\mP(2).
\]
Multiplying the equations in \eqref{eqn:prenew} and rearranging we have
\[
\lambda^\prime+\mu^\prime=1, \qquad \lambda^\prime:=\frac{\lambda_1^2}{\lambda_2^2},
\quad
\mu^\prime:=\frac{-4 \varepsilon}{\lambda_2^2}.
\]
Here we have a new solution $(\lambda^\prime,\mu^\prime)$ to the $S$-unit equation
\eqref{eqn:sunit}.
Moreover,
\[
\ord_{\mP}(\mu^\prime)=0, \qquad 
\ord_\mP(\lambda^\prime)=
2 \ord_{\mP}(\lambda)-4 \ord_\mP(2) > \ord_{\mP}(\lambda),
\]
where the last inequality follows from~\eqref{eqn:tlm}.
This contradicts the maximality
of $\ord_\mP(\lambda)$ and completes the proof.


\section{Examples and Comparisons}\label{sec:examples}
We wrote a short \texttt{Magma}
implementation of the criterion 
of Theorem~\ref{thm:compCrit}. We would of course like
to compare  this criterion with the criterion of Theorem~\ref{thm:main}. 
For $n\ge 3$ let $\cF_n$ be the set of totally real fields of degree $n$,
discriminant $\le 10^6$,
in which $2$ totally ramifies. By the theorem of Odlyzko
(quoted in \cite[Proposition 2.3]{Takeuchi}),
\[
29.009^n \cdot \exp{(-8.3185)} < 10^6.
\] 
It follows that $3 \le n \le 6$. We were able
to find the complete sets $\cF_3$, $\cF_4$, $\cF_5$, $\cF_6$
in the John Jones \texttt{Number Field Database} \cite{JR}. It turns
out that $\cF_6$ is empty, so we focus on degree $3$, $4$, $5$.
We define the following sets.
\begin{itemize}
\item Let $\cG_n$ be the set of $K \in \cF_n$
such that $h_K^+$ is odd.
\item Let $\cH_n$ be the set of $K \in \cF_n$ 
such that $[\mP] \in \Cl^+(K)$ has order
divisible by $h_{K,2}^+$ and $h_K$ is odd.
\item Let $\cI_n$ be the set of $K \in \cF_n$
such that $[\mP] \in \Cl^+(K)$ has order
divisible by $h_{K,2}^+$.
\item Let $\cJ_n$ be the set of $K \in \cF_n$
such that every element of $U$ is a square (where
we take $V=\OO_K^\times$).
\end{itemize}
Of course $\cG_n \subseteq \cH_n \subseteq \cI_n$.
We know the asymptotic FLT holds for any $K$ belonging
to $\cI_n$ or $\cJ_n$, thanks to Theorems~\ref{thm:main}
and~\ref{thm:compCrit}. We computed these sets using our \texttt{Magma}
implementation. The results are summarised in the table.

\begin{table}[h]
\begin{tabular}{||c|c|c|c|c|c||}
\hline\hline
$n$ & $\#\cF_n$ & $\#\cG_n$ & $\#\cH_n$ & $\#\cI_n$ & $\#\cJ_n$\\ 
\hline\hline
$3$ & $8600$ & $3488$ & $3488$ & $3488$ & $7653$\\
\hline
$4$ & $1243$ & $1$ & $428$ & $446$ & $1039$\\
\hline
$5$ & $23$ & $13$ & $13$ & $13$ & $22$\\
\hline\hline
\end{tabular}
\end{table}
We make the following observations.
\begin{enumerate}
\item[(I)] As $2$ is totally ramified in all these fields,
$2\OO_K=\mP^{n}$. If $n=3$, $5$, then the order of $[\mP]$
in $\Cl^+(K)$ is odd. It follows that
$\cG_n=\cH_n=\cI_n$.
\item[(II)] For $n=3$, $4$, $5$, we found that
$\cI_n \subset \cJ_n$. An explanation for this
is given by Lemma~\ref{lem:explain} below.
\item[(III)] Note
that $\#\cG_4=1$. In other words, there is only one totally real quartic
field $K$ of discriminant $\le 10^6$ for which $h_K^+$ is
odd. This field is $K=\Q(\zeta_{16})^+$. This observation is explained by
Theorem~\ref{thm:explain}.
\item[(IV)] The fields in $\cH_3$, $\cH_4$, $\cH_5$ are precisely
the totally real fields with discriminant $\le 10^6$
that satisfy conditions (a), (b), (c) of Theorem~\ref{thm:main2}.
For such $K \in \cH_3 \cup \cH_4 \cup \cH_5$ we know that asymptotic FLT holds
over any $K^\prime$ that is a $2$-extension of $K$ unramified away from 
the unique prime above $2$; in particular, when $K^\prime=K \cdot \Q(\zeta_{2^r})^+$.
\end{enumerate}

\subsection*{A Variant where $2$ is not totally ramified}
We give an example to show how the proof of Theorem~\ref{thm:compCrit}
can still be useful in establishing asymptotic FLT, even
if $2$ does not totally ramify in the field. Let
$K$ be the number field generated by a root of
$x^5+x^4-12x^3-21x^2+x+5$. This is totally real
and has degree $5$. Moreover, $2$ is inert in this
field, and we let $\mP=2\OO_K$. Taking $V=\OO_K^\times$,
we checked that every element of $U$ is a square. The
proof of Theorem~\ref{thm:compCrit} shows that if $(\lambda,\mu)$
is a solution to the $S$-unit equation, then
after applying a suitable element of $\sS_3$,
$0 \le \ord_\mP(\lambda) \le 4$ and $\ord_\mP(\mu)=0$.
Thus we may write $\lambda=2^r \lambda^\prime$ where $0 \le r \le 4$,
and $\lambda^\prime$, $\mu$ are both units. Instead of solving
the $S$-unit equation \eqref{eqn:sunit} to apply Theorem~\ref{thm:FS},
we merely have to solve the (easier) unit equations 
$2^r \lambda^\prime+\mu=1$ for $0 \le r \le 4$.
Using \texttt{Magma}'s inbuilt unit equation solver
we find no solutions for $r=0$, $2$, $3$, $4$, and
precisely one solution for $r=1$, which $\lambda^\prime=1$, $\mu=-1$.
Hence the only solutions to the $S$-unit equation \eqref{eqn:sunit}
are the irrelevant ones. Now applying part (B) of Theorem~\ref{thm:FS}
allows us to deduce asymptotic FLT for $K$.

\subsection*{Another Computational Example}
Let $\{f_n\}$ be the sequence of polynomials in \eqref{eqn:fn}.
It is easy to see that $f_n$ is monic of degree $n$ and
belongs to $\Z[x]$. Let $\alpha$ be any root of $f_n$
and $K_n=\Q(\alpha)$. 
\begin{lem}\label{lem:Kn}
$K_n$ is a totally real field of degree $n$ in which
$2$ totally ramifies.
\end{lem}
\begin{proof}
Write $K=K_n$.
Let $L=K(\sqrt{-7})$. Let $\beta \in L$ be related to $\alpha$ by
\[
\beta=\frac{\alpha+\sqrt{-7}}{\alpha-\sqrt{-7}}, \qquad
\alpha=\sqrt{-7} \cdot \frac{(\beta+1)}{(\beta-1)}.
\]
Let $\pi_1=(1+\sqrt{-7})/2$ and $\pi_2=(1-\sqrt{-7})/2$. 
These are the primes above $2$ in $\Q(\sqrt{-7})$.

From $f_n(\alpha)=0$ it follows that $\beta^n=\pi_2/\pi_1$, so that $x^n - \pi_2/\pi_1$ is $\pi_2$-Eisenstein polynomial. We conclude
$\pi_1$, $\pi_2$ are totally ramified in the degree $n$
extension $L/\Q(\sqrt{-7})$, and $[L:\Q]=2n$.
Moreover,  $\pi_2/\pi_1$ has complex absolute value 1 so
$\sigma(\beta)$ lies on the unit
circle for all embeddings $\sigma: L \hookrightarrow \C$.
Note that the M\"{o}bius transformation
$z \mapsto \sqrt{-7} \cdot (z+1)/(z-1)$ 
maps $\beta$ to $\alpha$ and
transforms
the unit circle into the real line. Thus $\sigma(\alpha) \in \R$
for all the embeddings of $L$, so $K$ is totally real.
Moreover $L=K(\sqrt{-7})$ and so $K$ has degree $n$. 
Finally $2$ is not ramified in $\Q(\sqrt{-7})$,
and thus the primes of $K$ above $2$ are not
ramified in $L/K$. This allows us to deduce
that $2$ is totally ramified in $K$.
\end{proof}
We ran our \texttt{Magma} implementation of
the criterion of Theorem~\ref{thm:compCrit} for these fields
with $2 \le n\le 32$.  Here we took $V$ to be a subgroup of
odd index in $\OO_K^\times$. We found that every element of $U$ 
is a square for precisely the following values of $n$:
$1 \le n \le 6$, $8 \le n \le 14$, $15 \le n \le 20$, $23 \le n \le 27$,
$n=29$, $31$, $32$. By Theorem~\ref{thm:compCrit} asymptotic FLT
holds over $K_n$ for these values. We note in passing that
\texttt{Magma} can compute the full unit group of $K_n$
for $n \le 18$, but appears not to be able to do this
(unconditionally) for larger values of~$n$.

\subsection*{A comparison of theorems}
The following lemma explains observation (II) above.
 Moreover, it does show that the assumptions
of Theorem~\ref{thm:compCrit} are weaker than those of 
Theorem~\ref{thm:main}. Indeed, Theorem~\ref{thm:main} is 
theoretically useful, but Theorem~\ref{thm:compCrit}
is more powerful in practice.
\begin{lem}\label{lem:explain}
Let $K$ be a number field with exactly one prime $\mP$
above $2$. Suppose the order of $[\mP] \in \Cl_K^+$ is
divisible by $h_{K,2}^+$. 
Let $V$, $U$ be as in Theorem~\ref{thm:compCrit}. Then
every element of $U$ is a square.
\end{lem}
\begin{proof}
Let $L$ be the narrow Hilbert $2$-class field of $K$. 
By Proposition~\ref{prop:equiv},
$L$ has exactly one
prime above $\mP$. 
Let $u \in U$. We want to prove that $u$ is a square in $K$.
We will in fact first show that it is a square in $L$ and then
deduce that it is a square in $K$.

Note that $(1+\sqrt{u})/2$ satisfies the polynomial 
$f=X^2-X+(1-u)/4 \in \OO_K[x]$. Moreover,
the discriminant of $f$ is $u \in \OO_K^\times$.
It follows in particular that $K(\sqrt{u})/K$
is an extension that is unramified at the finite places,
and so $K(\sqrt{u}) \subseteq L$ by definition of $L$.
Suppose $K(\sqrt{u}) \ne K$,
and write $M=K(\sqrt{u})$. 
Then $f$ is the minimal polynomial of $\alpha:=(1+\sqrt{u})/2 \in \OO_M$. 
We will apply the Dedekind--Kummer theorem to show that $\mP$
splits in $M$. 
To do this, we first show
that $\OO_M=\OO_K[\alpha]$. Indeed, by \cite[Theorem 1.2.30]{Cohen}
\[
\frac{\OO_M}{\OO_K[\alpha]} \cong \frac{\OO_K}{\ga_1}
\oplus \cdots \oplus \frac{\OO_K}{\ga_r}
\]
where $\ga_i$ are ideals of $\OO_K$ satisfying $\ga_1 \mid \ga_2 \mid \cdots
\mid \ga_r$. Moreover, by \cite[page 79]{Cohen},
\[
\Disc(f) \cdot \OO_K= \Delta_{M/K} \cdot \ga_1^2 \cdots \ga_r^2,
\]
where $\Delta_{M/K} \subseteq \OO_K$ is the relative discriminant ideal for $M/K$. 
But $\Disc(f)=u \in \OO_K^\times$.
It follows that the $\ga_i=\OO_K$ and so $\OO_M=\OO_K[\alpha]$
as desired.
Now $f \equiv X(X-1) \pmod{\mP}$. It follows from
the Dedekind--Kummer theorem \cite[Proposition 2.3.9]{Cohen} 
that $\mP$ splits in $K(\sqrt{u})$. Thus
there cannot be exactly one prime above $\mP$ in $L \supseteq M$.
This contradiction shows that $K(\sqrt{u})=K$.
\end{proof}

\subsection*{An Explanation for (III)}
The following is a mild generalization of \cite[Th\'{e}or\`{e}me 5]{Kraus}.
It provides an explanation for observation (III) above.
\begin{thm}\label{thm:explain} Let
$K$ be a totally real field of degree $2^n$ for  $n \geq 1$.
Suppose $2$ totally ramifies in $K$, and $h_K^+$ is odd.
Then $K=\Q(\zeta_{2^{n+2}})^+$.
\end{thm}
\begin{proof} 
We first show that every totally positive unit in $K$ is a square. 
Let $d=2^n$ for the degree of $K$. 
Write $V=\OO_K^\times$, $V^2$ for the subgroup of squares in $V$,
and $V^+$ for the subgroup of totally positive units.
By \cite[Corollary 3.2.4]{Cohen},
\[
\frac{h^+_K}{h_K}=\frac{2^d}{[V:V^+]}.
\]
As $h^+_K$ is odd, we have $h^+_K=h_K$ and so $[V:V^+]=2^d$.
However, $V^2\subseteq V^+$ and by Dirichlet's unit theorem
$[V:V^2]=2^d$. We conclude that $V^+=V^2$. In other words,
every totally positive unit of $K$ is a square.

We now need some notation.
For $r \geq 3$, write
\[
\beta_r=\zeta_{2^r}+\zeta_{2^r}^{-1},\qquad
\gamma_r=\beta_r+2, \qquad
L_r=\Q(\beta_r)=\Q(\zeta_{2^r})^+.
\]
Recall that the unique prime above $2$ in $\Q(\zeta_{2^r})$ is generated
by $(1-\zeta_{2^r})$ and therefore also
by its Galois conjugate $(1+\zeta_{2^r})$. Thus the unique 
prime above $2$ in $L_r$
is generated by $\gamma_r=(1+\zeta_{2^r})(1+\zeta_{2^r}^{-1})$.
We note that
\begin{equation} \label{eqn:gamma}
\beta_{r+1}^2=\gamma_r, \qquad (\gamma_r \OO_{L_r})^{2^{r-2}} =2 \OO_{L_r}.
\end{equation}
It follows from these that $\gamma_r:=\beta_r+2$ is totally positive
and $\gamma_r^{2^{r-2}}/2$ is a unit. 

\smallskip

\textbf{Claim:} $\beta_r \in K$ for $3 \le r \le n+2$.

\smallskip

Note that our claim implies the theorem, for applying the claim
with $r=n+2$ yields $\Q(\zeta_{2^{n+2}})^+ \subseteq K$, and as both fields
have degree $2^n$, so they must be equal.

\smallskip

Let $\mP$ be the unique
prime of $K$ above~$2$. Then $\mP^{2^n}=2\OO_K$. As $h_K^+$ is odd, we see that
$\mP=\alpha \OO_K$, where $\alpha \in \OO_K$ is totally positive. Then
$2/\alpha^{2^n}$ is a totally positive unit.
In follows that $2/\alpha^{2^n}=s^2$ where $s$ is a unit of $K$.
We prove the claim by induction.
Suppose $n \ge 1$ (if $n=0$ then there is nothing to prove).
Thus $\sqrt{2}=s\cdot \alpha^{2^{n-1}} \in K$.
But $\beta_3=\zeta_8+\zeta_8^{-1}=\sqrt{2}$. This establishes
our claim for $r=3$.
For the inductive step, suppose $3 \le r \le n+1$ and $\beta_r \in K$.
Therefore $L_r \subseteq K$ and $\gamma_r \in K$. From \eqref{eqn:gamma}
\[
(\gamma_r \OO_K)^{2^{r-2}} =2 \OO_{K}=\mP^{2^n}=\left(\mP^{2^{n+2-r}}\right)^{2^{r-2}}.
\]
Therefore
\[
\gamma_r \OO_K=\mP^{2^{n+2-r}}=\alpha^{2^{n+2-r}} \OO_K.
\]
Thus $\gamma_r/\alpha^{2^{n+2-r}}$ is a unit.
It is totally positive, as $\alpha$ and $\gamma_r$ are totally positive,
and so must be the square of a unit, $s_r^2$.
Hence $\beta_{r+1}=\sqrt{\gamma_r}=s_r \cdot \alpha^{2^{n+1-r}}
\in K$, establishing the claim and completing the proof.
\end{proof}

\section{Proof of Theorem~\ref{thm:KrausconjII}}\label{sec:KrausconjII}

All the work done in this paper so far was under the assumption that the
fields~$K$ considered have a unique prime above~$2$. In this setting our
results were obtained by a careful study of the solutions to $S$-unit equation
with the help of class field theory.   

Theorem~\ref{thm:KrausconjII} concerns a family of quadratic
 fields having two primes
above $2$. We still study $S$-unit equation using class field theory, but that
alone appears insufficient to yield a complete proof. 
In particular, we will need the following proposition, whose
proof makes use of the theory of linear forms in logarithms and Diophantine
approximation. We postpone its proof to
Section~\ref{sec:dio} so as not to disrupt the flow of the argument.

\begin{prop}\label{prop:dio}
Let $\tau=3+2\sqrt{2}$. The only solutions to the equation
\begin{equation}\label{eqn:prell}
2^{s_1}+\eta \cdot 2^{s_2} = \frac{\tau^k-\tau^{-k}}{2\sqrt{2}}, \qquad
s_1,\, s_2,\, k \ge 0, \quad s_1 \ge s_2, \quad \eta=\pm 1,
\end{equation}
are $k=0$, $\eta=-1$, and $s_1=s_2$, or 
$(k,\eta,s_1,s_2)=(1,1,0,0)$, $(1,-1,2,1)$, $(2, 1, 3, 2)$, $(2, -1, 4, 2)$.
\end{prop}

We note that Theorem~\ref{thm:KrausconjII} follows immediately from Theorem~\ref{thm:FS} and the following lemma, so
the rest of this section is devoted to its proof.

\begin{lem}\label{lem:KrausConjII}
Let $\ell \equiv 1 \pmod{24}$ be  prime, and $K=\Q(\sqrt{\ell})$.
Write $\mP_1$, $\mP_2$ for the two primes of $K$ above $2$ and let
$S=\{\mP_1,\mP_2\}$. If $\ell>73$ then
the solutions to the $S$-unit equation
\begin{equation}\label{eqn:LambdaS}
\lambda+\mu=1, \qquad \lambda,\; \mu \in \OO_S^\times
\end{equation}
satisfy $\max\{\lvert \ord_{\mP_1}(\lambda) \rvert, \lvert \ord_{\mP_1}(\mu) \rvert\}
=1$ or $\max\{\lvert \ord_{\mP_2}(\lambda) \rvert, \lvert \ord_{\mP_2}(\mu) \rvert\}
=1$. 
If $\ell=73$, the same conclusion holds
 with the exception of the $\sS_3$-orbit of 
\begin{equation}\label{eqn:orbit}
\lambda=\frac{-23+3\sqrt{73}}{2}, \qquad \mu=\frac{25+3\sqrt{73}}{2}
\end{equation}
\end{lem}

From now on, the notation in this section 
will be that of Lemma~\ref{lem:KrausConjII}.

\begin{lem}\label{lem:ray}
The ray class numbers $h_{\mP_1^2}$, $h_{\mP_2^2}$ are odd.
\end{lem}
\begin{proof}
Let $\mP$ be either $\mP_1$ or $\mP_2$.
By \eqref{eqn:genus}, the class number $h_K$ of $K$ is odd.
Note that $(\OO_K/\mP^2)^\times \cong (\Z/4\Z)^\times$ is generated
by the image of $-1 \in \OO_K^\times$. Thus the natural
map
\[
\rho : \OO_K^\times \rightarrow (\OO_K/\mP^2)^\times
\]
is surjective. 
The exact sequence in  \cite[Proposition 3.2.3]{Cohen} 
tells us that $h_{\mP^2}=h_K$.
\end{proof}

Recall that the solutions $(1/2,1/2)$, $(-1,2)$ and
$(2,-1)$ to \eqref{eqn:LambdaS} are called
irrelevant, and the other solutions are called relevant. We would
like to understand the relevant solutions.
For the following, see \cite[Lemma 6.4]{FS} and its proof.

\begin{lem}\label{lem:param}
Up to the action of $\sS_3$, every relevant solution $(\lambda,\mu)$
of \eqref{eqn:LambdaS}
has the form
\begin{equation}\label{eqn:parasol}
\lambda=\frac{\eta_1 \cdot 2^{r_1}-\eta_2 \cdot 2^{r_2}+1 +v \sqrt{\ell}}{2},
\qquad 
\mu=\frac{\eta_2 \cdot 2^{r_2}-\eta_1 \cdot 2^{r_1}+1 -v \sqrt{\ell}}{2}
\end{equation}
where
\begin{equation}\label{eqn:paracond1}
\eta_1=\pm 1, \qquad \eta_2=\pm 1, \qquad r_1 \ge r_2 \ge 0, \qquad v \in \Z, \qquad v \ne 0
\end{equation}
are related by
\begin{gather}\label{eqn:paracond2}
(\eta_1 \cdot 2^{r_1}-\eta_2 \cdot 2^{r_2}+1)^2
-\eta_1 \cdot 2^{r_1+2}
=
\ell v^2 ,
\\
\label{eqn:paracond3}
(\eta_2 \cdot 2^{r_2}-\eta_1 \cdot 2^{r_1}+1)^2
-\eta_2 \cdot 2^{r_2+2}
=\ell v^2 .
\end{gather}
Moreover,
\begin{equation}\label{eqn:norm}
\Norm_{K/\Q}(\lambda)=\eta_1 \cdot 2^{r_1}, \qquad \Norm_{K/\Q}(\mu)=\eta_2 \cdot 2^{r_2}.
\end{equation}
\end{lem}
\begin{lem}\label{lem:r1large}
Let $(\lambda,\mu)$, $\eta_i$, $r_i$ be as in Lemma~\ref{lem:param}
with $r_1 \le 5$. Then $\ell=73$ and $(\lambda,\mu)$
is given by \eqref{eqn:orbit}.
\end{lem}
\begin{proof}
This is a straightforward computation as $r_2 \le r_1$.
\end{proof}

We shall henceforth suppose that $r_1 \ge 6$.
\begin{lem}\label{lem:parity}
Let $(\lambda,\mu)$, $\eta_i$, $r_i$ be as in Lemma~\ref{lem:param}.
Then, for $i=1$, $2$, we have
\[
\eta_i=-1 \iff r_i \equiv 1 \pmod{2}.
\]
\end{lem}
\begin{proof}
Let $a=\eta_1 \cdot 2^{r_1}$ and $b=\eta_2 \cdot 2^{r_2}$. 
Equation~\eqref{eqn:paracond2}
becomes
\[
(a-b+1)^2-4a=\ell v^2.
\]
We consider this modulo $3$. Since $\ell \equiv 1 \pmod{3}$
we infer 
\begin{equation}\label{eqn:simple}
(a-b+1)^2-a \; \equiv \; \text{$0$ or $1 \pmod{3}$}.
\end{equation}
However $(a,b) \equiv (\pm 1,\pm 1) \pmod{3}$. Of these four possibilities
for $(a,b)$ modulo $3$, the only one that satisfies \eqref{eqn:simple}
is $(a,b) \equiv (1,1) \mod{3}$. Thus $2^{r_i} \equiv \eta_i \pmod{3}$
for $i=1$, $2$. This gives the lemma.
\end{proof}
\begin{lem}\label{lem:posval}
Let $(\lambda,\mu)$, $\eta_i$, $r_i$ be as in Lemma~\ref{lem:param}.
Then $r_2>0$.
\end{lem}
\begin{proof}
Suppose $r_2=0$. By Lemma~\ref{lem:parity} we have $\eta_2=1$.
Now \eqref{eqn:paracond2} becomes
\[
2^{2r_1}-\eta_1 \cdot 2^{r_1+2} = \ell v^2.
\]
The $2$-adic valuation of the left-hand side is $r_1+2$, and this must
be even in view of the right-hand side.  
Moreover, by Lemma~\ref{lem:parity} we have $\eta_1=1$. 
Removing a factor of $2^{r_1+2}$ from both sides gives 
\[
2^{r_1-2}-1=\ell w^2, \qquad w \in \Z.
\]
This is impossible modulo $4$.
\end{proof}
Since the residue field of $\mP_i$ is $\F_2$ and $\lambda+\mu=1$
we have $\mP_i$ divides $\lambda$ or $\mu$ but not both. 
In particular $\max\{\ord_{\mP_i}(\lambda),\ord_{\mP_i}(\mu)\}\ge 1$
for $i=1$, $2$. As we would like to prove Lemma~\ref{lem:KrausConjII}
we shall suppose that $\max\{\ord_{\mP_i}(\lambda),\ord_{\mP_i}(\mu)\}\ge 2$
for $i=1$, $2$.
From \eqref{eqn:norm}, if $\mP_1$, $\mP_2$
both divide $\lambda$ then $r_2=0$, and if both divide $\mu$ then $r_1=0$,
which contradict Lemma~\ref{lem:posval}. Hence, after possibly swapping
$\mP_1$, $\mP_2$ we have,
\begin{equation}\label{eqn:vals}
\begin{cases}
\ord_{\mP_2}(\lambda)=\ord_{\mP_1}(\mu)=0,\quad
\ord_{\mP_1}(\lambda)=r_1,\quad
\ord_{\mP_2}(\mu)=r_2,\\
r_1 \ge r_2, \quad 
r_1 \ge 6, \quad r_2 \ge 2.
\end{cases}
\end{equation}
\begin{lem}\label{lem:elem}
Let $\eta=\pm 1$.
The only solutions to the equation $a^2-b^2= \eta \cdot 2^k$ in positive
odd integers $a$, $b$
are 
\[
a=2^{k-2}+\eta, \qquad
b=2^{k-2}-\eta
\]
 with $k \ge 3$.
\end{lem}
\begin{proof}
Observe that $a^2-b^2 \equiv 0 \pmod{8}$ and so $k \ge 3$.
It is sufficient to prove the lemma for $\eta=1$.
Then $(a+b)(a-b)=2^k$,
and so
\[
a+b=2^{s}, \qquad a-b=2^{t}, \qquad 1 \le t <s, \qquad s+t=k.
\]
Then $b=2^{s-1}-2^{t-1}$ and as $b$ is odd we have $t=1$, and so $s=k-1$.
\end{proof}

\begin{lem}\label{lem:sgn}
$\eta_1=\eta_2=-1$.
\end{lem}
\begin{proof}
Suppose $\eta_1=1$. Thus $r_1=2s$ by Lemma~\ref{lem:parity}.
From \eqref{eqn:paracond2} we have
\[
(2^{2s}+2^{s+1}+1-\eta_2 2^{r_2})(2^{2s}-2^{s+1}+1-\eta_2 2^{r_2})=\ell v^2.
\]
The two factors are coprime. Moreover, as $2s=r_1 \ge r_2$ we see that the
first factor is positive, and so the second must be positive. Hence
\[
\text{Case I:} \quad 2^{2s}+2^{s+1}+1-\eta_2 2^{r_2}=x^2, \qquad 2^{2s}-2^{s+1}+1-\eta_2 2^{r_2}=\ell y^2
\]
or
\[
\text{Case II:} \quad 2^{2s}+2^{s+1}+1-\eta_2 2^{r_2}=\ell y^2, \qquad 2^{2s}-2^{s+1}+1-\eta_2 2^{r_2}=x^2
\]
for some positive integers $x$, $y$. 

Let
\[
\eta_3=
\begin{cases}
1 & \text{Case I}\\
-1 & \text{Case II}.
\end{cases}
\]
Then we may rewrite our equations as
\begin{equation}\label{eqn:cases}
(2^s+\eta_3)^2-x^2=\eta_2 2^{r_2}, \qquad 
2^{2s}-\eta_3 2^{s+1}+1-\eta_2 2^{r_2}=\ell y^2.
\end{equation}
We apply Lemma~\ref{lem:elem} to the first equation in \eqref{eqn:cases}. 
This gives us
\[
2^s+\eta_3=2^{r_2-2}+\eta_2,
\]
and $r_2 \ge 3$. Recall that $r_1 \ge 6$. Thus $s \ge 3$.
Hence $\eta_3=\eta_2$ and $s=r_2-2$. Substituting into the second equation
in \eqref{eqn:cases} we obtain
\[
2^{2r_2-4}-3\eta_2 \cdot 2^{r_2-1}+1=\ell y^2.
\]
Thus $\ell y^2 \equiv 2 \pmod{3}$, which is impossible. This completes 
the proof that $\eta_1=-1$. 

It remains to show that $\eta_2=-1$. Thus suppose $\eta_2=1$, hence
$r_2$ is even and we write $r_2=2s$. Now
we use \eqref{eqn:paracond3}, which we can rewrite as
\[
(2^{2s}+2^{s+1}+1+2^{r_1})(2^{2s}-2^{s+1}+1+2^{r_1})=\ell v^2.
\] 
the factors are positive and coprime. Now the proof is exactly as before.
\end{proof}

From Lemma~\ref{lem:sgn} and Lemma~\ref{lem:parity} we know
now that $r_1$, $r_2$ are odd. We shall write $r_i=2s_i+1$.
We can now make \eqref{eqn:norm} more precise:
\begin{equation}\label{eqn:norm2}
\Norm(\lambda)=-2^{2s_1+1}, \qquad \Norm(\mu)=-2^{2s_2+1}.
\end{equation}

We shall denote Galois conjugation in $K$ by $x \mapsto \overline{x}$.
\begin{lem}\label{lem:totpos}
$\lambda \overline{\mu}$ is totally positive (i.e. positive in both embeddings).
\end{lem}
\begin{proof}
Let $\sigma_1$, $\sigma_2 : K \hookrightarrow \R$ be the two embeddings.
Let $\lambda_i=\sigma_i(\lambda)$ and $\mu_i=\sigma_i(\mu)$.
We are required to show that $\lambda_1 \mu_2>0$ and $\lambda_2 \mu_1>0$.
From \eqref{eqn:norm2}
\[
\lambda_1 \lambda_2<0, \qquad \mu_1\mu_2<0.
\]
Moreover as $\lambda+\mu=1$ we have $\lambda_1+\mu_1=1$
and $\lambda_2+\mu_2=1$. Thus $\lambda_1$, $\mu_1$ cannot both be negative,
and $\lambda_2$, $\mu_2$ cannot both be negative. Examining all the possible
signs we find that $\lambda_1 \mu_2$ and $\lambda_2 \mu_1$ are both positive.
\end{proof}

\begin{lem}
$\lambda \overline{\mu}=\varepsilon^2$ for some $\varepsilon \in \OO_K$.
\end{lem}
\begin{proof}
We shall study the ramification for $K(\sqrt{\lambda \overline{\mu}})/K$.
Note that $\lambda \overline{\mu}$ is totally positive by Lemma~\ref{lem:totpos},
thus $K(\sqrt{\lambda \overline{\mu}})/K$ is unramified at the infinite places.
Moreover, it is unramified at all finite $\mathfrak{q} \nmid 2$
as $\lambda \overline{\mu} \in \OO_S^\times$. 
To study ramification at $\mP_1$, $\mP_2$. Note
that $K_{\mP_1}=K_{\mP_2}=\Q_2$. There are seven
 quadratic extensions of $\Q_2$ obtained by adjoining
the square-root of one of $5$, $3$, $7$, $2$, $6$, $10$, $14$,
and these have discriminants $1$, $2^2$, $2^2$, $2^3$, $2^3$, $2^3$, $2^3$
respectively.
As $\mP_2=\overline{\mP_1}$, we see from the valuations in
\eqref{eqn:vals} that
\[
\lambda \overline{\mu}=(1-\mu)\cdot (1-\overline{\lambda}) \equiv 1 \pmod{\mP_2^{r_2}}.
\]
Thus
$\mP_2$
does not ramify in $K(\sqrt{\lambda \overline{\mu}})$. It remains
to measure the ramification at $\mP_1$. However,
$\ord_{\mP_1}(\lambda \overline{\mu})=r_1+r_2$. From Lemma~\ref{lem:sgn},
and Lemma~\ref{lem:parity}, we have $r_1 \equiv r_2 \equiv 1 \pmod{2}$,
and so $2 \mid (r_1+r_2)$.  
We deduce that the discriminant of $K(\sqrt{\lambda \overline{\mu}})/K$
divides $\mP_1^2$. As the ray class number for the modulus $\mP_1^2$
is odd (Lemma~\ref{lem:ray}), we deduce that $K(\sqrt{\lambda \overline{\mu}})=K$.
\end{proof}

We now complete the proof of Lemma~\ref{lem:KrausConjII}.
From \eqref{eqn:norm2}
\[
1=(\lambda+\mu)(\overline{\lambda}+\overline{\mu})=-2^{2s_1+1}-2^{2s_2+1}+\lambda \overline{\mu}+
\overline{\lambda} \mu.
\]
Thus
\begin{equation}\label{eqn:epstr}
\varepsilon^2+\overline{\varepsilon}^2=2^{2s_1+1}+2^{2s_2+1}+1.
\end{equation}
However,
\[
(\varepsilon \overline{\varepsilon})^2=\Norm(\varepsilon^2)=\Norm(\lambda \mu)=2^{2s_1+2s_2+2}.
\]
Thus
\begin{equation}\label{eqn:epsnr}
\varepsilon \overline{\varepsilon} =\eta \cdot 2^{s_1+s_2+1}, \qquad \eta=\pm 1.
\end{equation}
Recall $K=\Q(\sqrt{\ell})$ with $\ell \equiv 1 \pmod{24}$.
As $\varepsilon \in \OO_K$ we may write 
\[
\varepsilon=\frac{w_1+w_2 \sqrt{\ell}}{2}
\]
where $w_1$, $w_2 \in \Z$, and $w_1 \equiv w_2 \pmod{2}$.
Now \eqref{eqn:epstr} and \eqref{eqn:epsnr} can be rewritten as
\begin{equation}\label{eqn:pair}
\frac{w_1^2+\ell w_2^2}{2}=2^{2s_1+1}+2^{2s_2+1}+1, \qquad
\frac{w_1^2-\ell w_2^2}{2}=\eta \cdot 2^{s_1+s_2+2}.
\end{equation}
Hence
\begin{equation}\label{eqn:prepell}
w_1^2=2^{2s_1+1}+\eta \cdot 2^{s_1+s_2+2}+2^{2s_2+1}+1
=2 \cdot (2^{s_1}+\eta \cdot 2^{s_2})^2+1,
\end{equation}
and 
\begin{equation}\label{eqn:lw2}
\ell w_2^2=2^{2s_1+1}+2^{2s_2+1}+1 - \eta \cdot 2^{s_1+s_2+2}.
\end{equation}
We rewrite \eqref{eqn:prepell} as
\[
w_1^2 - 2 \cdot (2^{s_1}+\eta \cdot 2^{s_2})^2 = 1.
\]
It follows that $\lvert w_1 \rvert+(2^{s_1}+\eta \cdot 2^{s_2}) \sqrt{2}$
is a unit in $\Z[\sqrt{2}]$. The units of this ring have 
the form $\pm (1+\sqrt{2})^u$. However, as $1+\sqrt{2}$
has norm $-1$, we deduce that
\begin{equation}\label{eqn:pell}
\lvert w_1 \rvert +(2^{s_1}+\eta \cdot 2^{s_2}) \sqrt{2} =\tau^k,
\end{equation}
where $\tau=(1+\sqrt{2})^2=3+2\sqrt{2}$ and $k$ is a non-negative integer. Hence
\[
2^{s_1}+\eta \cdot 2^{s_2} = \frac{\tau^k-\tau^{-k}}{2\sqrt{2}}.
\]
We now apply Proposition~\ref{prop:dio} to deduce that
$k=0$, $\eta=-1$ and $s_1=s_2$ (the other solutions in
the proposition lead
either to the solution \eqref{eqn:orbit} or
to contradictions with \eqref{eqn:lw2} using $\ell \equiv 1 \pmod{24}$). 
From \eqref{eqn:lw2} we obtain
\[
\ell w_2^2=2^{2s_1+3}+1=8 \times 4^{s_1}+1.
\]
The right-hand side is divisible by $3$, and so $3 \mid w_2$. Thus the
right-hand side must be divisible by $9$. But the right-hand is $0$, $6$ or $3$
modulo $9$
according to whether $s_1 \equiv 0$, $1$ or $2 \pmod{3}$. Hence $3 \mid s_1$.
Write $s_1=3t$. Then
\[
\ell w_2^2=(2^{2t+1})^3+1=(2^{2t+1}+1)(2^{4t+2}-2^{2t+1}+1).
\]
The greatest common divisor of the two factors on the right-hand side is $3$.
Hence either
$2^{2t+1}+1=3x^2$,  or $2^{4t+2}-2^{2t+1}+1=3x^2$, for some integer $x$. However,
both equations are impossible modulo $4$, as long as $t \ge 1$. If $t=0$
then $\ell w_2^2=9$ which is also impossible.
This completes the proof of Lemma~\ref{lem:KrausConjII} and therefore
of Theorem~\ref{thm:KrausconjII}.

\section{A Diophantine Problem: Proof of Proposition~\ref{prop:dio}}\label{sec:dio}
In Section~\ref{sec:KrausconjII} we used class field
theory to reduce the proof of Theorem~\ref{thm:KrausconjII}
to the Diophantine problem in Proposition~\ref{prop:dio}.
We will now give a proof of that proposition.
Throughout this section $(k,\eta,s_1,s_2)$ will be a solution to \eqref{eqn:prell}.
The conclusion of Proposition~\ref{prop:dio} is clear when $k=0$. Therefore assume that $k \ge 1$. 

\begin{lem}\label{lem:twolog}
$s_2 \le \ord_2(k)+1$.
\end{lem}
\begin{proof}
We work in  $\Z[\sqrt{2}]$. We claim that
\[
\ord_{\sqrt{2}} (\tau^{2^a}-1)=2a+3, \qquad \text{for all $a \ge 1$}.
\]
This is true for $a=1$, and claim is easily established by induction
using the identity
\[
\tau^{2^{a+1}}-1=2\cdot \left(\tau^{2^a}-1\right)+\left(\tau^{2^a}-1\right)^2.
\]

Write $k=2^b \cdot k_0$ where $k_0$ is odd.
From \eqref{eqn:prell} we deduce that
\[
(\tau^{2^{b+1}})^{k_0}=\tau^{2k} \equiv 1 \pmod{\sqrt{2}^{2s_2+3}}.
\]
But, using the standard number field generalization of the
Euler totient function,
\[
\# \left(\Z[\sqrt{2}]/\sqrt{2}^{2s_2+3} \right)^\times=
\Norm\left(\sqrt{2}^{2s_2+3}\right)-
\Norm\left(\sqrt{2}^{2s_2+2}\right)=2^{2s_2+2},
\]
and in particular is coprime to $k_0$.  Thus
\[
\tau^{2^{b+1}}\equiv 1 \pmod{\sqrt{2}^{2s_2+3}}.
\]
We deduce that 
\[
2s_2+3 \le \ord_{\sqrt{2}} \left(\tau^{2^{b+1}}-1 \right)=2(b+1)+3,
\]
so $\ord_{2}(k)=b \ge s_2-1$ completing the proof.
\end{proof}

\begin{lem}\label{lem:k1000}
The only solutions to \eqref{eqn:prell} with $k \le 10^3$ are 
as given in Proposition~\ref{prop:dio}.
\end{lem}
\begin{proof}
Let $P_k=(\tau^k-\tau^{-k})/2\sqrt{2} \in \Z$. 
For each value $1 \le k \le 10^3$, 
Lemma~\ref{lem:twolog} gives us $s_2 \le \ord_2(k)+1$. 
In view of \eqref{eqn:prell}, we 
need only check for the possible values of $k$, $s_2$, and $\eta=\pm 1$,
if $P_k-\eta \cdot 2^{s_2}$ is a power of $2$.
We wrote a \texttt{Magma} script which did this
and the only solutions we found are
$(k,\eta,s_1,s_2)=(1,1,0,0)$, $(1,-1,2,1)$, $(2, 1, 3, 2)$, $(2, -1, 4, 2)$.
\end{proof}

We will therefore henceforth assume that $k>1000$. This makes
the forthcoming inequalities easier to deal with.
%
Next we apply the theory of linear forms in logarithms to obtain
a bound on $k$. We know that $2^{s_2} \le 2k$ by Lemma~\ref{lem:twolog}.
From \eqref{eqn:prell} we obtain 
\[
\left\lvert  \sqrt{2}^{2s_1+3} -\tau^k \right\rvert
\le 2^{s_2+1} \cdot \sqrt{2}+\tau^{-k}
\le 4 \sqrt{2} \cdot k+\tau^{-k} < 6k.
\]
Hence
\begin{equation}\label{eqn:prelin}
\left\lvert 
\frac{\sqrt{2}^{2s_1+3}}{\tau^k} \, - \, 1  
\right\rvert \; < \; 
\frac{6k}{\tau^k}\, .
\end{equation}
Let
\begin{equation}\label{eqn:Lam}
\Lambda=(2s_1+3) \cdot \log(\sqrt{2}) \, - \, k \cdot \log(\tau) \, .
\end{equation}
Using the elementary inequality $\lvert \log(1+x) \rvert \le 2 x$
for $\lvert x \rvert \le 1/2$, where we take $x=\exp(\Lambda)-1$,
we obtain
\begin{equation}\label{eqn:bd}
\left\lvert 
\Lambda
\right\rvert \; < \; \frac{12k}{\tau^k}.
\end{equation}
Hence
\begin{equation}\label{eqn:2sp3}
2s_1+3 \; <  \;
\left(\frac{\log(\tau)}{\log(\sqrt{2})} + \frac{12}{\log(\sqrt{2}) \cdot \tau^{1000}} \right) \cdot k \; < \;
6  k\, .
\end{equation}
We now apply the theorem of Baker and W\"{u}stholz for linear forms in logarithms
\cite[page 225]{Smart}, where in the notation of that theorem we take:
\[
\alpha_1=\sqrt{2}, \quad \alpha_2=\tau, \quad n=d=2, \quad b_1=2s_1+3, \quad b_2=-k, \quad B=\max\{2s_1+3,k\}.
\]
We find in the notation of that theorem (see also \cite[page 22]{Smart}) that
\[
h_m(\alpha_1)=\frac{1}{2}, \quad h_m(\alpha_2)=\frac{\log{(3+2\sqrt{2})}}{2}.
\]
The theorem gives
\begin{equation}\label{eqn:BW}
\log \lvert \Lambda \rvert > -C \cdot \log{B},
\end{equation}
where
\[
C=18 (n+1)! \cdot n^{n+1} \cdot (32 d)^{n+2} \cdot \log(2nd) \cdot h_m(\alpha_1) \cdot h_m(\alpha_2)
< 1.33 \cdot 10^{10}.
\]
From \eqref{eqn:2sp3} we have
\[
B < 6 k.
\]
Thus from \eqref{eqn:bd} and \eqref{eqn:BW} we obtain 
\[
\log{k}+\log{12}-k \cdot \log{\tau} > -C \cdot \log{k}-C\cdot \log{6}.
\]
Thus
\[
k \; < \; \frac{C+1}{\log{\tau}} \cdot \log{k}+\frac{C \cdot \log{6}+\log{12}}{\log{\tau}}
\; < \; 
a+b \log{k}
\]
where $a=
1.36 \times 10^{10}$ and $b=7.55 \times 10^9$.
Now Lemma B.1 of \cite[Appendix B]{Smart} gives
\begin{equation}\label{eqn:kbound}
k \; < \; 2(a+b\log{b}) \; < \; 3.8 \times 10^{11}.
\end{equation}

From \eqref{eqn:Lam} and \eqref{eqn:bd} we find
\begin{equation}\label{eqn:approx}
\left\lvert
\frac{2s_1+3}{k} 
\, - \,
\frac{\log{\tau}}{\log{\sqrt{2}}}
\right\rvert
\; < \; 
\frac{12}{\log{\sqrt{2}}} \cdot \frac{1}{\tau^k} \, .
\end{equation}
We computed using \texttt{PARI/GP} the first $30$ terms of the
continued fraction expansion
of $\log{\tau}/\log{\sqrt{2}}$, and found that the $30$-th convergent
is $p/q$ where
\[
p=1815871259660093, \qquad
q=357018312787640 \approx 3.57 \times 10^{14}.
\]
Then
\begin{equation}
\left\lvert
\frac{p}{q} 
\, - \,
\frac{\log{\tau}}{\log{\sqrt{2}}}
\right\rvert
\; < \; 
\frac{1}{q^2}.
\end{equation}
Therefore
\[
\left\lvert \frac{p}{q} \, - \, \frac{2s_1+3}{k} \right\rvert \; < \;  
\frac{1}{q^2} \, +\, \frac{12}{\log{\sqrt{2}}} \cdot \frac{1}{\tau^k} \, .
\]
If $p/q=(2s_1+3)/k$ then $q \mid k$ contradicting \eqref{eqn:kbound}.
Thus $p/q \ne (2s_1+3)/k$ and so $\lvert p/q-(2s_1+3)/k \rvert \ge 1/qk$.
Thus
\[
\frac{1}{qk} \; <  \; \frac{1}{q^2} \, +\, \frac{12}{\log{\sqrt{2}}} \cdot \frac{1}{\tau^k} \, .
\]
From \eqref{eqn:kbound}, $1/(2k)>1/q$ and so, 
\[
\frac{1}{2q k} \; < \; \frac{12}{\log{\sqrt{2}}} \cdot \frac{1}{\tau^k}
\]
and thus
\[
k  \; < \; \frac{1}{\log{\tau}} \cdot \log{k}+\frac{1}{\log{\tau}} \cdot \log\left(\frac{24q}{\log{\sqrt{2}}} \right).
\]
Applying Lemma B.1 of \cite[Appendix B]{Smart} now gives $k<73$.
In view of Lemma~\ref{lem:k1000}, this
completes the proof of Proposition~\ref{prop:dio}.


\begin{thebibliography}{99}
\bibitem{BF}
J.-F.\ Biasse and C.\ Fieker, 
\emph{Improved techniques for computing the ideal class group and a system of
fundamental units in number fields},
ANTS X—Proceedings of the
Tenth Algorithmic Number Theory Symposium, 113--133, Open Book Ser., 
\textbf{1}, Math.
Sci.\ Publ., Berkeley, CA, 2013.

\bibitem{BST}
M.\ Bhargava, A.\ Shankar, J.\ Tsimerman, 
\emph{On the Davenport-Heilbronn theorems and second order terms},  
Invent.\ Math.\ \textbf{193} (2013), no.\ 2, 439--499. 


\bibitem{magma}
W.\ Bosma, J.\ Cannon and C.\ Playoust,
\emph{The Magma algebra system I: the user language},
 J.\ Symb.\ Comp. \textbf{24} (1997), 235--265. (See also \url{http://magma.maths.usyd.edu.au/magma/})

\bibitem{browkin} J.\ Browkin, 
\emph{The abc-conjecture for Algebraic Numbers}, 
Acta Math.\ Sinica, English series {\bf 22} (2006), 211--222.


\bibitem{Buchmann}
J.\ Buchmann, 
\emph{A subexponential algorithm for the determination of class groups and
regulators of algebraic number fields}, 
S\'{e}minaire de Th\'{e}orie des Nombres, Paris
1988--1989, 27--41, Progr. Math. \textbf{91}, Birkh\"{a}user Boston, 
Boston, MA, 1990.




\bibitem{Cohen}
H.\ Cohen,
\emph{Advanced topics in computational number theory},
GTM \textbf{138}, Springer Verlag, 2000.

\bibitem{CDO}
H.\ Cohen, F.\ Diaz y Diaz and M. Olivier, 
\emph{Subexponential algorithms for class group and unit computations}, 
J.\ Symb.\ Comp. \textbf{24} (1997), 433--441.

\bibitem{Cohn}
H.\ Cohn, 
\emph{The density of abelian cubic fields},
Proc.\ Amer.\ Math.\ Soc.\ \textbf{5}, (1954) 476--477. 


\bibitem{FS} 
N.\ Freitas and S.\ Siksek, 
\emph{The asymptotic Fermat's last theorem for five-sixths of real quadratic fields},
Compos.\ Math.\ \textbf{151} (2015), no.\ 8, 1395--1415. 

\bibitem{FS2} 
N.\ Freitas and S.\ Siksek, 
\emph{Fermat's Last Theorem over some small real quadratic fields},
Algebra \& Number Theory {\bf 9} (2015), no. 4, 875--895.



\bibitem{Hanrot} G.\ Hanrot, 
\emph{Solving Thue equations without the full unit group}, 
Math.\ Comp.\ \textbf{69} (2000), no.\ 229, 395--405. 

\bibitem{HSV}
W.\ Ho, A.\ Shankar and I.\ Varma,
\emph{Odd degree number fields with odd class number},
Duke Math.\ J.\ \textbf{167} (2018), 995--1047.

\bibitem{Iwasawa} K.\ Iwasawa, 
\emph{A note on class numbers of algebraic number fields},
Abh.\ Math.\ Sem.\ Univ.\ Hamburg \textbf{20} (1956), 257--258. 

\bibitem{JarvisMeekin}
F.\  Jarvis and P.\  Meekin, 
\emph{The Fermat equation over $\Q(\sqrt{2})$},
J.\ Number Theory \textbf{109} (2004), no. 1, 182--196.

\bibitem{Kraus}
A.\ Kraus,
\emph{Le th\'{e}or\`{e}me de Fermat sur certains corps de nombres totalement r\'{e}els},
Algebra \& Number Theory, to appear.

\bibitem{LS} 
H.\ W.\ Lenstra, P.\ Stevenhagen,
\emph{Class Field Theory and the First Case of Fermat’s Last Theorem},
in G.\ Cornell, J.\ H.\ Silverman and G.\ Stevens (eds) 
\emph{Modular Forms and Fermat’s Last Theorem},
Springer, New York, 1997.

\bibitem{Lemmermeyer}
F.\ Lemmermeyer, 
\emph{Class Field Towers}, preprint, 111 pages, September 7, 2010.

\bibitem{JR}
J.\ W.\ Jones and D.\ P.\ Roberts, 
\emph{A database of number fields},
LMS J.\ Comput.\ Math.\ \textbf{17} (2014), no.\ 1, 595--618. 
(See also \url{https://hobbes.la.asu.edu/NFDB/})

\bibitem{Marksaitis}
G.\ N.\ Mark\u{s}a\u{\i}tis, 
\emph{On $p$-extensions with one critical number},
 Izv.\ Akad.\ Nauk SSSR Ser.\ Mat.\ \textbf{27} (1963), 463--466.

\bibitem{PARI}
The PARI~Group, PARI/GP version {\tt 2.11.0}, Univ. Bordeaux, 2018,
\url{http://pari.math.u-bordeaux.fr/}.

\bibitem{Pop}
F.\ Pop, 
\emph{Embedding problems over large fields}, 
Ann.\ of Math.\ (2) \textbf{144} (1996), no.\ 1, 1--34. 

\bibitem{Haluk} M.\ H.\ \c{S}eng\"{u}n and S.\ Siksek,
\emph{On the asymptotic Fermat's Last Theorem over numbers fields}, 
Commentarii Mathematici Helvetici
\textbf{93} (2018), 359--375.

\bibitem{SerreBook} J.-P.\ Serre,
\emph{Abelian $\ell$-adic representations and elliptic curves}, 
Addison-Wesley Publ. Co., Reading, Mass., 1989.

\bibitem{SilvermanII} J.\ H.\ Silverman,
{\em Advanced Topics in the Arithmetic of Elliptic Curves},
GTM \textbf{151}, Springer, 1994.

\bibitem{Smart}
N.\ P.\ Smart,
\emph{The Algorithmic Resolution of Diophantine Equations},
Cambridge University Press, 1998.


\bibitem{Suzuki}
M.\ Suzuki,
\emph{Group Theory I},
Grundlehren der mathematischen Wissenschaften \textbf{247},
Springer, Berlin, 1982. 



\bibitem{Takeuchi}
K.\ Takeuchi, 
\emph{Arithmetic Fuchsian groups with signature $(1;e)$}, 
J.\ Math.\ Soc.\ Japan \textbf{35} (1983), no.\ 3, 381--407.

\bibitem{Thorne}
J.\ Thorne,
\emph{Elliptic curves over $\Q_\infty$ are modular},
to appear in JEMS.

\bibitem{Turcas}
G.\ Turcas,
\emph{On Fermat’s equation over some quadratic imaginary number fields}, 
Res.\ Number Theory {\bf 4} (2018), no. 2, Art. 24, 16 pp. 

\bibitem{Washington}
L.\ C.\ Washington, 
\emph{Introduction to Cyclotomic Fields}, 
Graduate Texts in Mathematics \textbf{83}, Springer,
1982.

\bibitem{Wiles} A.\ Wiles, 
\emph{Modular elliptic curves and Fermat's last theorem}, Ann. of Math. (2) \textbf{141} (1995), no. 3, 443--551.
\end{thebibliography}
\end{document}